\documentclass[12pt,oneside]{amsart}
\usepackage[cp1250]{inputenc}
\usepackage[T1]{fontenc}
\usepackage{color}

\begin{document}
\newtheorem{theorem}{Theorem}[section]
\newtheorem{lemma}[theorem]{Lemma}
\newtheorem{remark}[theorem]{Remark}
\newtheorem{example}[theorem]{Example}
\newtheorem{corollary}[theorem]{Corollary}
\newtheorem{proposition}[theorem]{Proposition}

\begin{center}
{\LARGE\textbf{On piecewise-regular and quasirational functions}\vspace*{5mm}}\\
{\large Marcin Bilski}
\end{center}\vspace*{5mm}
\textbf{Abstract. } We fully characterize piecewise-regular functions on nonsingular real algebraic varieties in
terms of continuous quasirational functions. We also obtain extension theorem for piecewise-regular functions
defined
on closed semialgebraic subsets of $\mathbb{R}^n.$   \vspace*{2mm}\vspace*{2mm}

\section{Introduction}\label{secintro}

Continuous rational functions in real algebraic geometry have attracted attention of several mathematicians (cf.,
for example, \cite{FHMM}, \cite{KoKuKu}, \cite{KoNo}, \cite{Ku2014}, \cite{KuKu2}, \cite{MoJP}, \cite{Zie} and the
references therein). Interest in these functions is partially motivated by their approximation properties (see
\cite{Ku2014}, \cite{KuKu2}, \cite{Zie}). However, in many cases continuous rational functions turn out to be too
rigid to approximate arbitrary continuous maps between real algebraic varieties (cf. \cite{Ku2014},
\cite{Ku2017}). For this reason larger classes of maps contained in the class of semi-algebraic ones have been
studied. An example is the class of piecewise-regular functions which are much more flexible and proved to be
useful to approximate continuous maps into uniformly rational real algebraic varieties  (cf. \cite{BiKu},
\cite{Bi}, \cite{Ku2017}).

The aim of the present paper is to fully characterize piecewise-regular functions in terms of continuous
quasirational functions which apparently generalize continuous rational functions in a more straightforward way.
The developed techniques allow us to obtain extension theorem for piecewise-regular functions which is of
independent interest.

Before stating our main results let us note that real algebraic varieties appearing in the present paper are
affine and are endowed with the Zariski topology and the Euclidean topology. The notion of a nonsingular point is
understood in the algebraic sense. Moreover, given a real algebraic variety $V$, all its points with local
dimension smaller than $\mathrm{dim}(V)$ are considered to be singular. Unless explicitly stated otherwise, all
topological notions relating to subsets of $\mathbb{R}^n$ refer to the Euclidean topology induced by the standard
norm on $\mathbb{R}^n.$ For more information regarding real algebraic varieties we refer to either \cite{BCR} or
\cite{ManB}.\vspace*{2mm}

The zero-locus of a function $f:K\rightarrow\mathbb{R}$ will be denoted by $\mathcal{Z}(f).$\vspace*{2mm}

\noindent\textbf{Definition.} Let $K$ be a subset of $\mathbb{R}^n.$ A function $f:K\rightarrow\mathbb{R}$ is said
to be \textit{regular} if there are $P, Q\in\mathbb{R}[x_1,\ldots,x_n]$ such that $\mathcal{Z}(Q)\cap K=\emptyset$
and $$f(x)=\frac{P(x)}{Q(x)}\mbox{ for }x\in K.$$ \textbf{Definition.} Let $S$ be a semialgebraic subset of
$\mathbb{R}^n.$ A function $f:S\rightarrow\mathbb{R}$ is said to be \textit{piecewise-regular} if it is continuous
and there is a finite family $\{S_j\}_{j\in J}$ of semialgebraic subsets of $\mathbb{R}^n$ with $\bigcup_{j\in
J}S_j=S$ such that $f|_{S_j}$ is regular for $j\in J.$\vspace*{2mm}

The original definition of a piecewise-regular function is formulated in \cite{Ku2017} in a more general context
in a bit different way. However, on semi-algebraic sets both definitions are equivalent (see Proposition
\ref{PiecewiseSemi} below).\vspace*{2mm}

\noindent\textbf{Definition.} A continuous function $b:\mathbb{R}^n\rightarrow\mathbb{R}$ is said to be
\textit{basic quasipolynomial} if there exists $v\in\mathbb{R}[x_1,\ldots,x_n]$ such that $$|b(x)|=|v(x)|\mbox{
for all }x\in\mathbb{R}^n.$$ A continuous function $h:\mathbb{R}^n\rightarrow\mathbb{R}$ is said to be \textit{
quasipolynomial} if $h$ is a finite sum of basic quasipolynomial functions. The class of all quasipolynomial
functions on $\mathbb{R}^n$ will be denoted by $QP(\mathbb{R}^n).$\vspace*{2mm}\\
\textbf{Definition.} Let $V$ be a Zariski locally closed subvariety of $\mathbb{R}^n.$ A function
$g:V\rightarrow\mathbb{R}$ is said to be \textit{quasirational} if there exist $w\in\mathbb{R}[x_1,\ldots,x_n],
h\in QP(\mathbb{R}^n)$ such that $\mathcal{Z}(w)\cap V$ is nowhere dense in $V$ in the Euclidean topology and
$$g(x)=\frac{h(x)}{w(x)}, \mbox{ for }x\in V\setminus\mathcal{Z}(w).$$
The class of all quasirational functions on $V$ will be denoted by $QR(V).$\vspace*{2mm}

Our main results are the following two theorems

\begin{theorem}\label{quasipiece} Let $V$ be a nonsingular Zariski locally closed subvariety of $\mathbb{R}^n.$ Then every
continuous quasirational function $V\rightarrow\mathbb{R}$ is a piecewise-regular function.
\end{theorem}

\begin{theorem}\label{mainpw} Let $S$ be a closed semialgebraic subset of $\mathbb{R}^n.$ Then for every
piecewise-regular function $f:S\rightarrow\mathbb{R},$ there exists a continuous quasirational function
$g:\mathbb{R}^n\rightarrow\mathbb{R},$ $$g(x)=\frac{h(x)}{w(x)}, \mbox{ for }x\in
\mathbb{R}^n\setminus\mathcal{Z}(w),$$ where $h\in QP(\mathbb{R}^n), w\in\mathbb{R}[x_1,\ldots,x_n]$ such that
$\mathrm{dim}(\mathcal{Z}(w)\cap S)<\mathrm{dim}(S)$ and $g|_S=f.$
\end{theorem}

Note the we cannot skip the assumption that $S$ is closed. For nonclosed subsets even continuous extensions to
$\mathbb{R}^n$ in general do not exist.

The first corollary of Theorems \ref{quasipiece}, \ref{mainpw} says that on nonsingular real algebraic sets the
classes of piecewise-ragular and continuous quasirational functions coincide.

\begin{corollary}\label{corollary1} Let $V$ be a nonsingular real algebraic subset of $\mathbb{R}^n$. Then for every function
$f:V\rightarrow\mathbb{R}$ the following conditions are equivalent:\vspace*{2mm}\\
(a) $f$ is a piecewise-regular function,\\
(b) $f$ is a continuous quasirational function
\end{corollary}
\textit{Proof of Corollary~\ref{corollary1}.} The implication $(b)\Rightarrow (a)$ is an immediate consequence of
Theorem~\ref{quasipiece}.

Let now $f$ be a piecewise-regular function on $V.$ By Theorem~\ref{mainpw}, $f$ can be extended to a
quasirational function $g$ equal to $\frac{h}{w}$ on $\mathbb{R}^n\setminus\mathcal{Z}(w),$ where
$\mathrm{dim}(\mathcal{Z}(w)\cap V)<\mathrm{dim}(V).$ Since all irreducible components of $V$ are of the same
dimension we have that $\mathcal{Z}(w)\cap V$ is nowhere dense in $V$ in the Euclidean topology. Hence $f\in
QR(V)$ as required.\qed\\
\\
\noindent\textbf{Remark.}{ Theorem \ref{quasipiece} and the implication $(b)\Rightarrow(a)$ in
Corollary~\ref{corollary1} are in general not true unless we assume that $V$ is nonsingular. Namely, for singular
$V$ the restriction $f|_{\mathrm{Sing}(V)}$ of a continuous rational function $f:V\rightarrow\mathbb{R}$ may be a
non-rational Nash function (see \cite{KoNo}). Then $f$ cannot be piecewise-regular.\vspace*{2mm}

The two following corollaries of Theorems~\ref{quasipiece}, \ref{mainpw} deal with the extension problems.
\begin{corollary}\label{extpiece} Let $S$ be a closed semialgebraic subset of $\mathbb{R}^n$ and $f:S\rightarrow\mathbb{R}$
a piecewise-regular function. Then there is a piecewise-regular function $F:\mathbb{R}^n\rightarrow\mathbb{R}$
such that $F|_S=f.$
\end{corollary}
\textit{Proof of Corollary \ref{extpiece}.} By Theorem \ref{mainpw}, $f$ can be extended to a continuous
quasirational function on $\mathbb{R}^n$, which by Theorem \ref{quasipiece} is piecewise-regular.\qed\vspace*{2mm}

The following result is also an immediate consequence of Theorems \ref{quasipiece}, \ref{mainpw}.
\begin{corollary}\label{extquasi} Let $V$ be a nonsingular real algebraic subset of $\mathbb{R}^n$  and $g:V\rightarrow\mathbb{R}$ a continuous quasirational function. Then there is
a continuous quasirational function $G:\mathbb{R}^n\rightarrow\mathbb{R}$ such that $G|_V=g.$
\end{corollary}

Again, Corollary \ref{extquasi} is not true unless we assume that $V$ is nonsingular. Indeed, if the continuous
rational function from Example 2 of \cite{KoNo} had a continuous quasirational extension to $\mathbb{R}^n,$ then
the extension would be a piecewise-regular function on $\mathbb{R}^n$ by Theorem \ref{quasipiece}. Consequently,
the original function would be piecewise-regular which is not possible (remark above).

\section{Auxiliary results I}\label{introduction}


For any subset $A$ of $\mathbb{R}^n,$ let $\overline{A}^z$ be the closure of $A$ in the Zariski topology of
$\mathbb{R}^n.$ Let us recall the notion of piecewise-regular map from \cite{Ku2017}.

A \textit{stratification} of a real algebraic variety $V$ is, by definition, a finite collection of pairwise
disjoint Zariski locally closed subvarieties whose union equals
$V.$\vspace*{0mm}\\
\\
\textbf{Definition.} Let $X$ be any subset of $\mathbb{R}^n,$ $f:X\rightarrow\mathbb{R}$ a continuous map and
$\mathcal{S}$ a stratification of $\mathbb{R}^n.$ The map $f$ is said to be \textit{piecewise
$\mathcal{S}$-regular} if for every stratum $S\in\mathcal{S}$ the restriction of $f$ to each connected component
of $X\cap S$ is a regular map (when $X\cap S$ is non-empty). Moreover, $f$ is said to be
\textit{piecewise-regular} if it is piecewise $\mathcal{T}$-regular for some stratification $\mathcal{T}$ of
$\mathbb{R}^n.$\vspace*{2mm}

It is clear that if $X$ is semialgebraic and $f$ is piecewise-regular in the sense of the definition above
(Section\ref{introduction}), then it is also piecewise-regular in the sense of the definition in Section
\ref{secintro}.

\begin{proposition}\label{PiecewiseSemi} Let $S$ be a semialgebraic subset of $\mathbb{R}^n$ and
$f:S\rightarrow\mathbb{R}$ a continuous function. Assume that there exists a finite family $\{S_{l}\}_{l\in L}$ of
semialgebraic subsets of $S$ such that $\bigcup_{l\in L}S_l=S$ and $f|_{S_l}$ is regular for each $l\in L.$ Then
$f$ is piecewise-regular in the sense of definition above (Section \ref{introduction}).
\end{proposition}

\textit{Proof.} Without loss of generality we assume that the elements of $\{S_l\}_{l\in L}$ are pairwise disjoint
(as $\{S_l\}_{l\in L}$ can be replaced by semialgebraic covering of $S$ by pairwise disjoint sets each of which is
contained in some $S_l$).

We need a stratification of $\overline{S}^z$ such that for each stratum $C,$ the restriction of $f$ to every
connected component of $C\cap S$ is regular. The existence of such a stratification will be proved by induction on
$\mathrm{dim}(S).$ If $\mathrm{dim}(S)=0,$ there is nothing to prove.

For every $S_l$ with $\mathrm{dim}(S_l)=\mathrm{dim}(S),$ we have $S_l=Y_l\cup Z_l,$ where $Y_l, Z_l$ are disjoint
semialgebraic sets such that $\mathrm{dim}(Y_l)<\mathrm{dim}(S)$ and $Z_l$ is open in~${S}.$ After replacing each
$S_l$ with $\mathrm{dim}(S_l)=\mathrm{dim}(S)$ by two sets $Y_l$ and $Z_l$ (and after reindexing  the elements of
$\{S_l\}_{l\in L}$) we may assume that every element of $\{S_l\}_{l\in L}$ of maximal dimension is open in $S.$

Put $$M:=\{l\in L: \mathrm{dim}(S_l)<\mathrm{dim}(S)\}\mbox{ and }Y:=\bigcup_{l\in M}S_l.$$ We define one of the
strata of $\overline{S}^z$ to be $\overline{S}^z\setminus\overline{Y}^z.$ To see that this is correct, fix a
connected component $A$ of $(\overline{S}^z\setminus\overline{Y}^z)\cap S.$ Put
$$N:=\{l\in L:\mathrm{dim}(S_l)=\mathrm{dim}(S)\}$$ and observe that $A\subset\bigcup_{l\in N } S_{l}.$ It follows
that $A$ is contained in some $S_l,$ as $\{S_l\}_{l\in N }$ consists of pairwise disjoint open subsets of $S.$
Hence $f$ is regular on $A,$ as required.

Put $W:=S\cap\overline{Y}^z.$ Clearly $Y\subset W,$ so $\overline{Y}^z\subset\overline{W}^z.$ Moreover,
$W\subset\overline{Y}^z$ so $\overline{W}^z=\overline{Y}^z.$ We have
$\mathrm{dim}(W)=\mathrm{dim}({Y})<\mathrm{dim}(S).$ Furthermore, $$W=\bigcup_{l\in L}(S_l\cap\overline{Y}^z)$$
and $f$ is regular on each $S_l\cap \overline{Y}^z.$ Hence, by induction hypothesis applied to $W$, there exists a
stratification $\mathcal{C}$ of $\overline{W}^z=\overline{Y}^z$ such that for every $C\in\mathcal{C},$ the
function $f$ is regular on each connected component of $C\cap W.$

Since $C\cap W=C\cap S,$ for $C\in\mathcal{C},$ the family
$$\{\overline{S}^z\setminus\overline{Y}^z\}\cup\mathcal{C}$$ is the stratification of $\overline{S}^z$ we look
for.\qed
\begin{corollary}\label{corpSemi} Let $S$ be a semialgebraic subset of $\mathbb{R}^n$ and
$f:S\rightarrow\mathbb{R}$ a continuous function. Assume that there exists a finite family $\{S_{l}\}_{l\in L}$ of
semialgebraic subsets of $S$ such that $\bigcup_{l\in L}S_l=S$ and $f|_{S_l}$ is piecewise-regular for each $l\in
L.$ Then $f$ is piecewise-regular.
\end{corollary}
\textit{Proof.} By definition of piecewise-regular function, each element of $\{S_l\}_{l\in L}$ can be subdivided
into a finite family of semialgebraic sets on which $f$ is regular. Therefore the assertion follows immediately by
Proposition \ref{PiecewiseSemi}.\qed
\\


\begin{lemma}\label{overset} Let $V$ be an affine nonsingular real algebraic variety of dimension $m.$
Let $Z_1,\ldots,Z_p, W$ be real algebraic subsets of $V$ such  that $\mathrm{dim}(Z_j)=k<m$ and $Z_j$ is
irreducible for $j=1,\ldots,p$ and $\mathrm{dim}(W)<m.$ Then there is a real algebraic subset $Y$ of $V,$ such
that each irreducible component of $Y$ is of dimension $k+1$ and
$$Z_j\subset Y \mbox{ and } Z_j\setminus\mathrm{Sing}(Y)\neq\emptyset, \mbox{ for }j=1,\ldots,p,$$
$$T\setminus W\neq\emptyset\mbox{ for every irreducible component $T$ of $Y$}.$$
\end{lemma}
\textit{Proof.} First consider the case $V=\mathbb{R}^m.$ Let
$\pi:\mathbb{R}^{m-1}\times\mathbb{R}\rightarrow\mathbb{R}^{m-1}$ be the natural projection. Let
$Z:=Z_1\cup\ldots\cup Z_p.$ We may assume (applying a linear change of variables in $\mathbb{R}^m$ if necessary)
that $\pi|_{Z\cup W}$ is a proper map. Then $\pi(Z_j)$ is a semialgebraic subset of $\mathbb{R}^{m-1}$ such that
$$k=\mathrm{dim}(\pi(Z_j))=\mathrm{dim}(\overline{\pi(Z_j)}^z), \mbox{ for }j=1,\ldots,p.$$
Note that the real algebraic set $\overline{\pi(Z_j)}^z$ is irreducible as $Z_j$ is irreducible.

Put $Y_j:=\overline{\pi(Z_j)}^z\times\mathbb{R}.$ Clearly, the set $Y_j$ is irreducible and
$$\mathrm{dim}(Y_j)=k+1\mbox{ and }Z_j\subset Y_j.$$
Define $Y:=Y_1\cup\ldots\cup Y_p.$ Then $Y=\overline{\pi(Z)}^z\times\mathbb{R}$ and
$$\mathrm{dim}(\pi(Z_j))=\mathrm{dim}(\overline{\pi(Z)}^z)>\mathrm{dim}(\mathrm{Sing}\overline{\pi(Z)}^z).$$
Consequently, $$\pi(Z_j)\setminus\mathrm{Sing}\overline{\pi(Z)}^z\neq\emptyset$$ and therefore
$$Z_j\setminus\mathrm{Sing}(\overline{\pi(Z)}^z\times\mathbb{R})=Z_j\setminus\mathrm{Sing}(Y)\neq\emptyset\mbox{ for }j=1,\ldots,p.$$

Finally, $T\setminus W\neq\emptyset$ for every irreducible component $T$ of $Y$ as all fibers of $\pi|_W$ are
finite whereas $\pi|_T$ has infinite fibers.

Now, let $V$ be any affine nonsingular real algebraic variety of dimension $m.$ Without loss of generality we
assume that $V$ is a real algebraic subset of $\mathbb{R}^n=\mathbb{R}^m\times\mathbb{R}^{n-m}.$ Let
$\varrho:\mathbb{R}^m\times\mathbb{R}^{n-m}\rightarrow\mathbb{R}^m$ denote the natural projection. Applying a
linear change of variables, we may assume that $\varrho|_V$ is a proper map with finite fibers and that
$\mathrm{Reg}(Z_j), \mathrm{Reg}(W),$ $j=1,\ldots,p,$ are not contained in the set of points $x\in V$ at which
$d_x(\varrho|_V)$ is singular.

Define $\tilde{Z}_j:=\overline{\varrho(Z_j)}^z, \tilde{W}:=\overline{\varrho(W)}^z.$ Then $\tilde{Z}_j,$
$\tilde{W}$ are real algebraic subsets of $\mathbb{R}^m,$ $\mathrm{dim}(\tilde{Z}_j)=\mathrm{dim}({Z}_j)=k$ and
$\mathrm{dim}(\tilde{W})=\mathrm{dim}(W),$ $\tilde{Z}_j$ irreducible, for $j=1,\ldots,p.$ We have already proved
that there exists an algebraic subset $\tilde{Y}$ of $\mathbb{R}^m$ such that each irreducible component of
$\tilde{Y}$ has dimension $k+1$ and
$$\tilde{Z}_j\subset \tilde{Y} \mbox{ and } \tilde{Z}_j\setminus\mathrm{Sing}(\tilde{Y})\neq\emptyset,\mbox{ for }j=1,\ldots,p,$$
and $$\tilde{T}\setminus\tilde{W}\neq\emptyset,$$ for every irreducible component $\tilde{T}$ of $\tilde{Y}.$

Clearly, the real algebraic set $(\tilde{Y}\times\mathbb{R}^{n-m})\cap V$ has dimension at most $k+1.$ Let us
define $Y$ to be the union of all $k+1$-dimensional irreducible components of
$(\tilde{Y}\times\mathbb{R}^{n-m})\cap V.$ We show that $Y$ enjoys the required properties.

First, we check  that $Z_j\subset Y,$ for $j=1,\ldots,p.$ Fix $j.$ Since $\mathrm{Reg}(Z_j)$ is not contained in
the set of points $x\in V$ at which $d_x(\varrho|_V)$ is singular, there exists an open subset $U$ of $V$ such
that $\varrho|_{U}:U\rightarrow\varrho(U)$ is a Nash diffeomorphism and $U\cap\mathrm{Reg}(Z_j)\neq\emptyset.$ In
particular, $$\mathrm{dim}(\varrho(U\cap \mathrm{Reg}(Z_j)))=k.$$ On the other hand, from
$\tilde{Z}_j\setminus\mathrm{Sing}(\tilde{Y})\neq\emptyset$ and the irreducibility of $\tilde{Z}_j$ it follows
that $\mathrm{dim}(\tilde{Z}_j\cap\mathrm{Sing}(\tilde{Y}))\leq k-1.$ Consequently,
$$\mathrm{dim}(\varrho(U\cap\mathrm{Reg}(Z_j))\cap\mathrm{Sing}(\tilde{Y}))\leq k-1,$$ and (in view of
$\varrho(Z_j)\subset\tilde{Y}$) there exists $$\hspace*{-43mm}(\sharp)\hspace*{50mm}b\in
\varrho(U\cap\mathrm{Reg}(Z_j))\cap\mathrm{Reg}(\tilde{Y})\neq\emptyset.
$$
For sufficiently small open neighborhood $E$ of $b$ in $\varrho(U),$ we have
$$ E\cap\varrho(U\cap Z_j)\subset\varrho(U\cap \mathrm{Reg}(Z_j))\mbox{ and }E\cap\tilde{Y}\subset\mathrm{Reg}(\tilde{Y}).$$
We may assume that $E\cap\tilde{Y}$ is connected. Then $M:=\varrho|_U^{-1}(E\cap\tilde{Y})$ is a connected Nash
manifold of dimension $k+1$ and
$$\varrho|_U^{-1}(b)\in M\cap Z_j=\varrho|_U^{-1}(E)\cap Z_j=\varrho|_U^{-1}(E)\cap\mathrm{Reg}(Z_j).$$
Hence, $M\cap Z_j$ has nonempty interior in $\mathrm{Reg}(Z_j).$ It follows that $Z_j\subset \overline{M}^z$ as
$Z_j$ is irreducible. Moreover, $\overline{M}^z$ is a $k+1$-dimensional irreducible (as $M$ connected Nash
manifold) subvariety of $(\tilde{Y}\times\mathbb{R}^{n-m})\cap V$ which in turn has dimension at most $k+1,$ so
$\overline{M}^z$ is an irreducible component of $Y.$ Thus $Z_j\subset Y,$ as required.

Let us check that, for the fixed $j,$ we have $Z_j\setminus\mathrm{Sing}(Y)\neq\emptyset.$ By $(\sharp),$ the
point $b$ is regular in $\tilde{Y}.$ Clearly, $\varrho|_U^{-1}(b)\in U$ is a regular point of $V,$ and $V\cap U$
is a graph of a Nash map over $\varrho(U).$ Therefore, $\varrho|_U^{-1}(b)$ is a regular point of
$(\tilde{Y}\times\mathbb{R}^{n-m})\cap V,$ whose dimension in a neighborhood of $\varrho|_U^{-1}(b)$ equals $k+1.$
Hence, $\varrho|_U^{-1}(b)$ is a regular point of $Y.$ Again, by $(\sharp),$ we have $\varrho|_U^{-1}(b)\in Z_j,$
so $Z_j\setminus\mathrm{Sing}(Y)\neq\emptyset.$

Finally, suppose that some irreducible  component $T$ of $Y$ is contained in $W.$ Then
$\overline{\varrho(T)}^z\subset\tilde{W}.$ But $\overline{\varrho(T)}^z$ is an algebraic set of dimension $k+1$
contained in $\tilde{Y},$ so it is an irreducible component of $\tilde{Y}.$ A contradiction. The proof is
complete.\qed \vspace*{2mm}


\begin{lemma}\label{closeopen} Let $M\subset N$ be Nash manifolds of dimensions $k, k+1,$ respectively. Let $B\subset
N$ be an open semialgebraic subset such that $\overline{B}^N\cap M$ has  nonempty interior in $M.$ Let $p,
q:N\rightarrow \mathbb{R}$ be Nash functions such that $\mathcal{Z}(q)$ is nowhere dense in $N$ and
$\frac{p}{q}|_{\overline{B}^N\setminus\mathcal{Z}(q)}$ can be continuously extended to $\overline{B}^N.$ Then
there is an open subset $U$ of $N$ such that $U\cap M\neq\emptyset$ and $\frac{p}{q}|_{U\setminus\mathcal{Z}(q)}$
can be continuously extended to $U.$
\end{lemma}
\textit{Proof.} We assume that $B\cap M=\emptyset,$ as otherwise we can take $U=B.$

By hypothesis, there is an open set $E$ in $N$ such that
$$\hspace*{-52mm}(i)\hspace*{52mm}\overline{B}^N\cap M\cap E=M\cap E\neq\emptyset.$$ The problem is local, so we may
assume that $N=\mathbb{R}^{k+1},$ $E=(-c,c)^{k+1}$ for some $c>0,$  and
$$\hspace*{-50mm}(ii)\hspace*{50mm}M\cap E=\{x\in E: x_{k+1}=0\}.$$

Denote $$E_<:=\{x\in E:x_{k+1}<0\}\mbox{ and }E_>:=\{x\in E:x_{k+1}>0\}.$$ Since $B$ is semialgebraic, we have,
after shrinking $E$ if necessary, that $E_<\subset B$ or $E_>\subset B.$ Indeed, using the cell decomposition of
$B$ and having in mind that the problem is local, we may assume that $B$ is of the form
$$B=\{(\tilde{x},x_{k+1})\in \tilde{B}\times\mathbb{R}: u(\tilde{x})<x_{k+1}<v(\tilde{x})\},$$ where $\tilde{B}$
is an open connected set in $\mathbb{R}^k$ and $u,v:\tilde{B}\rightarrow\mathbb{R}$ are Nash functions. By $(i),
(ii)$ and the fact that $B\cap M=\emptyset,$ we have $u=0$ or $v=0,$ which easily implies the assertion. Without
loss of generality assume that $E_<\subset B.$

We have $$q|_{E}=x_{k+1}^s\tilde{q}\mbox{ and } p|_{E}=x_{k+1}^r\tilde{p},$$ where $\tilde{p},
\tilde{q}:E\rightarrow \mathbb{R}$ are Nash functions with $$(M\cap
E)\setminus\mathcal{Z}(\tilde{q})\neq\emptyset\neq (M\cap E)\setminus\mathcal{Z}(\tilde{p}).$$ Since $E_{<}\subset
B,$ the hypothesis implies that $\frac{p}{q}|_{\overline{E_<}^N\setminus\mathcal{Z}(q)}$ can be continuously
extended to $\overline{E_<}^N.$ Therefore $r\geq s.$ Consequently,
$$\frac{p}{q}|_{E\setminus\mathcal{Z}(q)}=\frac{x_{k+1}^{r-s}\tilde{p}}{\tilde{q}}|_{E\setminus\mathcal{Z}(q)}.$$ Now the assertion follows immediately
by the fact that $(M\cap E)\setminus\mathcal{Z}(\tilde{q})\neq\emptyset.$\qed

\begin{lemma}\label{globalcont} Let $N$ be a nonsingular Zariski locally closed subvariety of $\mathbb{R}^n$ (with all
irreducible components of the same dimension). Let $M$ be an irreducible Zariski closed hypersurface in $N$ and
let $p, q\in\mathbb{R}[x_1,\ldots,x_n]$ with $\mathcal{Z}(q)\cap N$ nowhere dense in $N.$ Assume that there exists
an open subset $U$ of $N$ such that $U\cap\mathrm{Reg}(M)\neq\emptyset$ and
$\frac{p}{q}|_{U\setminus{\mathcal{Z}(q)}}$ has continuous extension to $U.$ Then there is an algebraic subset
$\Xi$ of $M$ with $\mathrm{dim}(\Xi)<\mathrm{dim}(M)$ such that for every $x\in M\setminus\Xi$ there exists an
open neighborhood $E$ in $N$ for which $\frac{p}{q}|_{E\setminus\mathcal{Z}(q)}$ has continuous extension to $E.$
\end{lemma}
\textit{Proof of Lemma \ref{globalcont}.} We may assume that $M\subset\mathcal{Z}(q)$ as otherwise
$\Xi:=\mathcal{Z}(q)\cap M$ satisfies the required properties.

Let $A$ be the set of points $z\in M$ which do not admit an open neighborhood $E_z$ in $N$ to which $\frac{p}{q}$
extends continuously. Then $A\subset M$ is semialgebraic. Hence it is sufficient to show that
$A\cap\mathrm{Reg}(M)$ is nowhere dense in $\mathrm{Reg}(M)$ and to define $\Xi:=\overline{A}^z.$

Fix any $y_0\in\mathrm{Reg}(M).$ We will check that every open neighborhood of $y_0$ in $\mathrm{Reg}(M)$ contains
a point which does not belong to $A$ which will complete the proof.

Fix any $x_0\in U\cap \mathrm{Reg}(M),$ where $U$ is an open subset of $N$ such that
$U\cap\mathrm{Reg}(M)\neq\emptyset$ and $\frac{p}{q}|_{U\setminus{\mathcal{Z}(q)}}$ has continuous extension to
$U.$ Let $N^{\mathbb{C}},$ $M^{\mathbb{C}}$ be the complexifications of $N, M$ in $\mathbb{C}^n,$ respectively.
Then $\mathrm{Reg}(M^{\mathbb{C}})$ is connected (as the complex variety $M^{\mathbb{C}}$ is irreducible which in
turn follows by irreducibility of the real variety $M$). Since $N$ is nonsingular we have
$N\cap\mathrm{Sing}(N^{\mathbb{C}})=\emptyset.$ Consequently, $\mathrm{Reg}(M^{\mathbb{C}})$ cannot be contained
in $\mathrm{Sing}(N^{\mathbb{C}}),$ so $\mathrm{Reg}(M^{\mathbb{C}})\setminus\mathrm{Sing}(N^{\mathbb{C}})$ is
connected and dense in $\mathrm{Reg}(M^{\mathbb{C}}).$

Let $\gamma:[0,1]\rightarrow \mathrm{Reg}(M^{\mathbb{C}})\setminus\mathrm{Sing}(N^{\mathbb{C}})$ be a continuous
function such that $\gamma(0)=x_0$ and $\gamma(1)=y_0.$ Let $U_t$ be an open neighborhood of $\gamma(t)$ in
$N^{\mathbb{C}}$ such that $M^{\mathbb{C}}\cap U_t$ is connected and has a minimal defining function $\chi_t$ on
$U_t,$ for $t\in [0,1].$

Now consider $p,q$ as elements from $\mathbb{C}[x_1,\ldots,x_n].$ Since $M\subset \mathcal{Z}(q),$ it holds
$M^{\mathbb{C}}\subset\mathcal{Z}(q).$ We have $$p|_{U_t}=\chi_t^{r_t}\tilde{p}_t\mbox{ and
}q|_{U_t}=\chi_t^{s_t}\tilde{q}_t,$$ where $\tilde{p}_t, \tilde{q}_t: U_t\rightarrow\mathbb{C}$ are holomorphic
functions nonvanishing identically on $M^{\mathbb{C}}\cap U_t.$

Note that $r_0\geq s_0,$ as $\frac{p}{q}|_{U\setminus\mathcal{Z}(q)}$ has continuous extension to $U.$ Define
$$I:=\{t\in [0,1]: r_t\geq s_t\}.$$
We have $I\neq\emptyset,$ as $0\in I.$ Moreover, $I$ is both open  and closed subset of $[0,1].$ This follows
immediately by the fact that for the open and dense subset $C_t:=(M^{\mathbb{C}}\cap
U_t)\setminus(\mathcal{Z}(\tilde{p}_t)\cup \mathcal{Z}(\tilde{q}_t) )$ of $M^{\mathbb{C}}\cap U_t$ the following
holds: if $r_t\geq s_t,$ then for every $a\in C_t$ the quotient $\frac{p}{q}$ is bounded in some neighborhood of
$a$ in $U_t$; if $r_t<s_t,$ then for every $a\in C_t$ the quotient $\frac{p}{q}$ is unbounded in every
neighborhood of $a$ in $U_t.$ Hence $I=[0,1],$ in particular, $r_1\geq s_1.$

Since $\tilde{q}_1$ does not vanish identically on $M^{\mathbb{C}}\cap U_1,$ we have that
$\mathcal{Z}(\tilde{q}_1)\cap \mathrm{Reg}(M)\cap U_1$ is nowhere dense in $\mathrm{Reg}(M)\cap U_1.$ It follows
that every open neighborhood of $y_0$ in $\mathrm{Reg}(M)$ contains a point $y_1\in \mathrm{Reg}(M)\cap U_1$ which
does not belong to $A.$\qed

\section{Proof of Theorem~\ref{quasipiece}}
Let $V$ be a nonsingular Zariski locally closed subvariety of $\mathbb{R}^n.$ Let $g:V\rightarrow\mathbb{R}$ be a
continuous quasirational function. According to Proposition \ref{PiecewiseSemi} and Corollary~\ref{corpSemi}, to
prove Theorem~\ref{quasipiece}, it is sufficient to show that there is a finite semialgebraic cover of $V$ such
that $g$ is (piecewise-)regular on each of its elements.

There are basic quasipolynomial functions $h_1,\ldots,h_t:\mathbb{R}^n\rightarrow\mathbb{R}$ and
$w\in\mathbb{R}[x_1,\ldots,x_n]$ such that $\mathcal{Z}(w)\cap V$ is nowhere dense in $V$ (in the Euclidean
topology) and
$$g(x)=\frac{h_1(x)+\cdots +h_t(x)}{w(x)} \mbox{ for }x\in V\setminus\mathcal{Z}(w).$$
Let $q_j\in\mathbb{R}[x_1,\ldots,x_n]$ such that $|h_j(x)|=|q_j(x)|, \mbox{ for }x\in\mathbb{R}^n, j=1,\ldots,t.$

Clearly, $V\setminus\mathcal{Z}(w)$ can be divided into a finite union $\bigcup_{a\in A}T_a$ of pairwise disjoint
semialgebraic sets such that either $h_j|_{T_a}=q_j|_{T_a}$ or $h_j|_{T_a}=-q_j|_{T_a}$ or $h_j|_{T_a}=0$ for
$j=1,\ldots,t$ and $a\in A.$ Therefore $g$ is regular on every $T_a.$

To complete the proof, it is sufficient to show the following\vspace*{2mm}\\
\textbf{Claim.}\textit{\label{claim} There exists a finite semialgebraic cover $\mathcal{S}$ of
$\mathcal{Z}(w)\cap V$ such that
$g$ is piecewise-regular on each $S\in\mathcal{S}.$}\vspace*{2mm}\\
\textit{Proof of Claim.} We use induction on the dimension $k$ of $\mathcal{Z}(w)\cap V.$ If $k=0,$ there is
nothing to prove. Assume that $k>0$ and that Claim is true if $\mathcal{Z}(w)$ intersects the domain of $g$ in
dimension $l<k.$

The proof will be divided into two steps. In Step 1 we will construct a semialgebraic set $S\subset
\mathcal{Z}(w)\cap V$ such that $g$ is piecewise-regular on $S$ and $$(\mathcal{Z}(w)\cap V)\setminus S\subset
H\subset\mathcal{Z}(w)\cap V ,$$ where $H$ is an algebraic subset of $V$ with $\mathrm{dim}(H)<k.$ To prove that
$g$ is piecewise-regular on $S$ we will use the result of \cite{KoNo} saying that on nonsingular real algebraic
varieties continuous rational functions are hereditarily rational.

In Step 2, we complete the proof of Claim using the induction hypothesis.\vspace*{2mm}\\
\textbf{Step 1.} Let $Z_1,\ldots, Z_p$ be all irreducible components of $\mathcal{Z}(w)\cap V$ of dimension $k.$
Put $Z:=Z_1\cup\ldots\cup Z_p$ and note that $\mathcal{Z}(w)\cap V=Z\cup\mathrm{Sing}(\mathcal{Z}(w)\cap V)$ (i.e.
lower dimensional irreducible components are, by definition, in the singular locus). The semialgebraic set $S$ we
look for will be contained in $Z.$

First, we reduce the problem to the case where (a large portion of) $Z$ is of codimension $1$ in some nonsingular
Zariski locally closed subvariety of $V$ on which $w$ does not vanish identically. More precisely, using
Lemma~\ref{overset} with $W=\emptyset$ we get a real algebraic subset $Y$ of $V$ such that each irreducible
component of $Y$ is of dimension $k+1$ and
$$Z_j\subset Y \mbox{ and } Z_j\setminus\mathrm{Sing}(Y)\neq\emptyset, \mbox{ for }j=1,\ldots,p.$$
Define $X:=Y\setminus(\mathrm{Sing}(Y)\cup \mathrm{Sing}(Z)).$ Then $X$ is a nonsingular  affine real algebraic
variety of dimension $k+1$ and $Z \cap X$ is a nonsingular hypersurface in $X$.

Next, we extract from $X$ an open dense semialgebraic subset $X^0$ such that the numerator of $g$ is a polynomial
function on each connected component of $X^0.$ Namely, for each connected component $C$ of $X$ define
$$\tilde{C}:=C\setminus\bigcup\{\mathcal{Z}(q_j): q_j|_C\neq 0, j\in\{1,\ldots,t\}\}.$$ Let $X^0$ be the union of all $\tilde{C}.$
Clearly, $X^0$ is an open dense (in the Euclidean topology) semialgebraic subset of $X.$ Let $\mathcal{A}$ be the
(finite) family of all connected components of $X^0$. Elements of $\mathcal{A}$ are open in $X.$ Moreover, for
each $A\in\mathcal{A}$ and $j=1,\ldots,t$ the polynomial $q_j$ has a constant sign on $A$ or ${q_j}|_A=0.$

Let
$${\mathcal{C}}:=\{\overline{A}^X:A\in\mathcal{A}\}.$$ By the density of $X^0$ in $X,$ we get
$\bigcup{\mathcal{C}}=X.$ For each $C\in\mathcal{C}, j=1,\ldots,p,$ let $$A_{C,j}:=\mathrm{int}_{Z_j\cap X}(
Z_j\cap C).$$ Since $\bigcup\mathcal{C}=X$ and each element of $\mathcal{C}$ is semialgebraic, we have
$$\hspace*{-36mm}(x)\hspace*{36mm}\bigcup\{A_{C,j}:C\in\mathcal{C}, j=1,\ldots,p\}=(Z\cap X)\setminus \Sigma,$$ where $\Sigma$ is a semialgebraic subset of
$Z\cap X,$ $\mathrm{dim}(\Sigma)\leq k-1.$

The set $S,$ we look for, will be contained in $T:=(Z\cap X)\setminus\Sigma.$ Before discussing  the details let
us check that $$\mathrm{dim}((\mathcal{Z}(w)\cap V)\setminus T)\leq k-1.$$ First, recall that
$$\mathcal{Z}(w)\cap V=Z\cup \mathrm{Sing}((\mathcal{Z}(w)\cap V)).$$ Thus it is sufficient to observe that
$\mathrm{dim}(Z\setminus T)\leq k-1.$ We have
$$Z\setminus T\subset \Sigma\cup(Z\setminus X)\subset\Sigma\cup (Z\cap \mathrm{Sing}(Y))\cup \mathrm{Sing}(Z),$$
proving the observation as $\mathrm{dim}(Z\cap \mathrm{Sing}(Y))\leq k-1$ and $\mathrm{dim}\mathrm{Sing}(Z)\leq
k-1.$

Our aim is to show that for every $C\in\mathcal{C}$ and $j=1,\ldots,p,$ there is a semialgebraic subset
$\Sigma_{C,j}$ of $A_{C,j}$, $\mathrm{dim}(\Sigma_{C,j})\leq k-1,$ such that $g|_{A_{C,j}\setminus\Sigma_{C,j}}$
is piecewise-regular. Once this is done, $(x)$ and Corollary \ref{corpSemi} imply piecewise-regularity of $g$ on
$$S=T\setminus\bigcup_{C\in\mathcal{C},j\in\{1,\ldots,p\}}\Sigma_{C,j}.$$ Then defining $H$ to be the closure of
$(\mathcal{Z}(w)\cap V)\setminus S $ in the Zariski topology of $V$ we complete Step 1.

Fix $C\in\mathcal{C}, j\in\{1,\ldots,p\}$ such that $A_{C,j}\neq\emptyset$ and recall that $C=\overline{B}^X$ for
some $B\in\mathcal{A}.$ Since each $q_i$ has constant sign on $B$ or $q_i|_B=0,$ there exists
$P\in\mathbb{R}[x_1,\ldots,x_n]$ such that
$g|_{B\setminus\mathcal{Z}(w)}=\frac{P}{w}|_{B\setminus\mathcal{Z}(w)}.$ Continuity of $g$ on $X$ and the fact
that $\overline{B\setminus\mathcal{Z}(w)}^X=C$ imply that  $\frac{P}{w}|_{B\setminus\mathcal{Z}(w)}$ extends
continuously to $g|_C$ on $C.$

We need to show that there is a Zariski open subset $X_{C,j}$ of $X$ such that
$\frac{P}{w}|_{X_{C,j}\setminus\mathcal{Z}(w)}$ has continuous extension to $X_{C,j}\rightarrow\mathbb{R}$ and
$$X_{C,j}\cap A_{C,j}=A_{C,j}\setminus\Sigma_{C,j},$$ where $\Sigma_{C,j}$ is a semialgebraic subset
of $A_{C,j},$ $\mathrm{dim}(\Sigma_{C,j})\leq k-1.$ Once such an $X_{C,j}$ exists, Step 1 will be complete.
Indeed, $X_{C,j}$ is a nonsingular real algebraic variety so, by \cite{KoNo}, the continuous extension of
$\frac{P}{w}|_{X_{C,j}\setminus\mathcal{Z}(w)}$ to $X_{C,j}\rightarrow\mathbb{R}$ is a hereditarily rational
function. In particular, it is a piecewise-regular function. Since $$A_{C,j}\setminus \Sigma_{C,j}\subset
X_{C,j}\cap C$$ and the extension of $\frac{P}{w}|_{C\setminus\mathcal{Z}(w)}$ to $C$ equals $g|_C,$ the function
$g|_{A_{C,j}\setminus\Sigma_{C,j}}$ is piecewise-regular, as required.

To obtain $X_{C,j}$ first apply Lemma \ref{closeopen} with $N=X$ and $M=Z_j\cap X.$ We get an open subset
$U_{C,j}$ of $X$ such that $U_{C,j}\cap (Z_j\cap X)\neq\emptyset$ and
$\frac{P}{w}|_{U_{C,j}\setminus\mathcal{Z}(w)}$ has continuous extension to $U_{C,j}\rightarrow\mathbb{R}.$ This
allows us to apply Lemma \ref{globalcont} with $N=X$ and $M=Z_j\cap X.$ We obtain an algebraic subset $\Xi_{C,j}$
of $Z_j\cap X,$ $\mathrm{dim}(\Xi_{C,j})\leq k-1,$ such that for every $x\in (Z_j\cap X)\setminus\Xi_{C,j},$ there
is an open neighborhood $E$ in $X$ such that $\frac{P}{w}|_{E\setminus\mathcal{Z}(w)}$ has continuous extension to
$E\rightarrow\mathbb{R}.$ Finally, define
$$X_{C,j}:=X\setminus(\bigcup_{l\neq j}Z_l\cup\mathrm{Sing}(\mathcal{Z}(w)\cap V)\cup\Xi_{C,j})$$ and note that
$X_{C,j}\cap\mathcal{Z}(w)\subset (Z_j\cap X)\setminus\Xi_{C,j}.$ Consequently, the function
$\frac{P}{w}|_{X_{C,j}\setminus\mathcal{Z}(w)}$ has continuous extension to $X_{C,j}\rightarrow\mathbb{R}.$ It
remains to check that $X_{C,j}\cap A_{C,j}$ has small complement in $A_{C,j}.$ Since $A_{C,j}\subset Z_j\cap X$
and

$$Z_j\cap X_{C,j}=Z_j\cap
X\setminus(\mathrm{Sing}(\mathcal{Z}(w)\cap V)\cup \Xi_{C,j} ),$$ for
$$\Sigma_{C,j}:=A_{C,j}\cap(\mathrm{Sing}(\mathcal{Z}(w)\cap V)\cup \Xi_{C,j} )$$ we have
$\mathrm{dim}(\Sigma_{C,j})\leq k-1$ and $A_{C,j}\cap X_{C,j}=A_{C,j}\setminus\Sigma_{C,j}.$ Step 1 is
complete.\vspace*{2mm}\\
\textbf{Step 2.} By Step 1, there is a semialgebraic set $S\subset \mathcal{Z}(w)\cap V$ such that $g$ is
piecewise-regular on $S$ and $$(\mathcal{Z}(w)\cap V)\setminus S\subset H\subset \mathcal{Z}(w)\cap V,$$ where $H$
is an algebraic subset of $V$
 with $\mathrm{dim}(H)=l\leq k-1.$ Hence, in view of Corollary~\ref{corpSemi}, to complete the proof of Claim, it is sufficient to prove that $g$
 is piecewise-regular on $H.$ To do this, we will use the following lemma with $Z=H.$


\begin{lemma}\label{lowdim} Let $V$ be a nonsingular Zariski locally closed subvariety of $\mathbb{R}^n.$
Let $w\in\mathbb{R}[x_1,\ldots,x_n]$  such that $\mathcal{Z}(w)\cap V$ is nowhere dense in $V.$ Let $Z$ be a real
algebraic subset of dimension $l$ of  $\mathcal{Z}(w)\cap V.$ Then there exist nonsingular Zariski locally closed
subsets $V_1,\ldots V_t$ of $V$ such that $Z\subset V_1\cup\ldots\cup V_t$ and $\mathcal{Z}(w)\cap V_i$ is nowhere
dense in $V_i$ and
$$\mathrm{dim}(\mathcal{Z}(w)\cap V_i)\leq l, \mbox{ for }i=1\,\ldots,t.$$
\end{lemma}
\textit{Proof of Lemma \ref{lowdim}. } We may assume that $V$ is irreducible (as a real algebraic variety) because
otherwise we consider every irreducible component separately.

The proof by induction on $l.$ Let $l=0.$ Then $Z=Z_1\cup\ldots\cup Z_p$ where each $Z_j$ is a one-element set.
Applying Lemma \ref{overset} with  $W=\mathcal{Z}(w)\cap V$ we obtain a real algebraic curve $Y$ in $V$ such that
$Z\subset Y\setminus\mathrm{Sing}(Y)$ and no irreducible component $T$ of $Y$ is contained in $\mathcal{Z}(w).$
Define $V_1=Y\setminus\mathrm{Sing}(Y).$ Clearly, being a finite set, $\mathcal{Z}(w)\cap V_1$ is of dimension at
most $l=0$ and is nowhere dense in $V_1.$

Let $l>0.$ Assume Lemma \ref{lowdim} for $l-1.$ Let $Z_1,\ldots, Z_p$ be all irreducible components of $Z$ of
dimension $l.$ We apply Lemma \ref{overset} with $W=\mathcal{Z}(w)\cap V$ to obtain a real algebraic subset $Y$ of
$V$ of dimension $l+1$ satisfying the assertions of this lemma. Then $V_1:=Y\setminus\mathrm{Sing}(Y)$ is a
nonsingular Zariski locally closed subset of $V$ such that $Z_j\setminus\mathrm{Sing}(Y)\subset V_1,$ for
$j=1,\ldots,p,$ and $\mathcal{Z}(w)\cap V_1$ is nowhere dense in $V_1$ and $\mathrm{dim}(\mathcal{Z}(w)\cap
V_1)\leq l.$

Let $Z_{(l-1)}$ be the union of all irreducible components of $Z$ of dimension at most $l-1.$ Note that
$\mathrm{dim}(Z_j\cap\mathrm{Sing}(Y))\leq l-1$ for $j=1,\ldots,p,$  and define
$$\tilde{Z}=Z_{(l-1)}\cup(\bigcup_{j=1}^pZ_j\cap\mathrm{Sing}(Y)).$$ Then $Z=(\bigcup_{j=1}^pZ_j\setminus\mathrm{Sing}(Y))\cup\tilde{Z}$
and $\tilde{Z}$ is a real algebraic subset of $\mathcal{Z}(w)\cap V$ of dimension at most $l-1.$

By the induction hypothesis Lemma \ref{lowdim} is true with $Z$ replaced by $\tilde{Z}.$ The family $V_2,\ldots,
V_t$ obtained by applying Lemma \ref{lowdim} for $\tilde{Z}$ together with $V_1$ defined above satisfies the
required properties.\qed\vspace*{2mm}

By Lemma \ref{lowdim}, there exist nonsingular Zariski locally closed subsets $V_1,\ldots V_t$ of $V$ such that
$H\subset V_1\cup\ldots\cup V_t$ and $\mathcal{Z}(w)\cap V_i$ is nowhere dense in $V_i$ and
$$\mathrm{dim}(\mathcal{Z}(w)\cap V_i)\leq k-1, \mbox{ for }i=1\,\ldots,t.$$
By the induction hypothesis Claim is true with $V$ replaced by $V_i,$ for $i=1,\ldots,t.$ This shows that $g$ is
piecewise-regular on every $V_i\cap H.$ Again, by Corolary \ref{corpSemi}, $g$ is piecewise-regular on $H$ which
proves Claim. \qed\vspace*{2mm}\\
Thus the proof of Theorem \ref{quasipiece} is also complete.

\section{Auxiliary results II }

\subsection{Preparation}

Let $\pi:\mathbb{R}^{n-1}\times\mathbb{R}\rightarrow \mathbb{R}^{n-1}$ be the natural projection.

\begin{lemma}\label{localsur}Let $p\in(\mathbb{R}[x_1,\ldots,x_{n-1}])[x_n]$ be a monic polynomial in $x_n$ of
degree $r$ in $x_n.$ Let $$w(x_1,\ldots,x_n)=\prod_{j=0}^r\frac{\partial^j p}{\partial x_n^j}(x_1,\ldots,x_n).$$
Then for every $a\in\mathcal{Z}(w)$ and every open ball $B$ in $\mathbb{R}^n$ with $a\in B,$ the following
conditions are satisfied:\vspace*{2mm}

\noindent(1)  there exist open sets $U, V$ in $\mathbb{R}^{n-1}, \mathbb{R},$ respectively, such that $a\in
U\times V\subset B$ and
$$\pi(\mathcal{Z}(w)\cap(U\times V))=U$$

\noindent(2) $\mathcal{Z}(w)\cap B$ is a semialgebraic set of dimension $n-1.$
\end{lemma}
\textit{Proof.} Fix $a\in\mathcal{Z}(w)$ and an open ball $B$ in $\mathbb{R}^{n-1}\times\mathbb{R}$ with $a\in B.$
Assume that (1) does not hold. Then  the condition $$\frac{\partial^i p}{\partial x_n^i}(a)=0\mbox{ for some
}i\in\{0,\ldots,r-1\}$$ implies that
$$\frac{\partial^{i+1} p}{\partial
x_n^{i+1}}(a)=0.$$ Indeed, if $\frac{\partial^{i+1} p}{\partial x_n^{i+1}}(a)\neq 0,$ then by the implicit
function theorem, the set $$\{x\in\mathbb{R}^n:\frac{\partial^i p}{\partial x_n^i}(x)=0\}\subset\mathcal{Z}(w)$$
contains the graph of a function $U\rightarrow V$ for some open sets $U, V$ in $\mathbb{R}^{n-1},\mathbb{R},$
respectively, with $a\in U\times V\subset B.$ This means that (1) holds, contrary to our assumption.

Now, note that $\frac{\partial^r p}{\partial x_n^r}(a)\neq 0$ because $p$ is a monic polynomial in $x_n$ of degree
$r.$ On the other hand, $\frac{\partial^i p}{\partial x_n^i}(a)=0\mbox{ for some }i\in\{0,\ldots,r-1\}$ because
$a\in\mathcal{Z}(w).$ The obtained contradiction proves (1).

As for (2), the inequality $\mathrm{dim}(\mathcal{Z}(w)\cap B)\leq n-1$ holds because $w$ is a nonzero polynomial.
The other inequality is an immediate consequence of (1).\qed\vspace*{3mm}

For $i=0,\ldots,n-1,$ let $\pi_i:\mathbb{R}^i\times\mathbb{R}\rightarrow\mathbb{R}^i$ be the natural projection.
For $d=1,\ldots,n-1,$ define $\tilde{\pi}_d:\mathbb{R}^n\rightarrow\mathbb{R}^d,$
$\tilde{\pi}_d:=\pi_d\circ\pi_{d-1}\circ\cdots\circ\pi_{n-1}.$ Put $\tilde{\pi}_n:=\mathrm{Id}_{\mathbb{R}^n}.$
For $p\in\mathbb{R}[x_1,\ldots,x_{i+1}],$ let $\mathrm{deg}_{x_{i+1}}p$ be the degree of $p$ in $x_{i+1}.$

\begin{proposition}\label{formulation}Let $\{S_j\}_{j\in J}$ be a finite family of semialgebraic subsets of
$\mathbb{R}^n.$ Then  after a linear change of coordinates in $\mathbb{R}^n,$ there are real algebraic subsets
$W_0,\ldots, W_n$ of $\mathbb{R}^n,$ $W_n=\mathbb{R}^n,$ such that {for} every $i=0,\ldots,n-1$ and $j\in J,$ the
following hold:\vspace*{2mm}
$$\hspace*{-29mm}(1_i)\hspace*{29mm}\mathrm{dim}(W_i)=i\mbox{ and }W_i\subset W_{i+1},\mbox{ and
}\mathrm{Sing}(W_{i+1})\subset W_{i}.$$
$$\hspace*{-54.5mm}(2_{i,j})\hspace*{54mm} \partial_{W_{i+1}}(W_{i+1}\cap S_j)\subset W_i.$$
$(3_{i})$ there exists a $q_i\in(\mathbb{R}[x_1,\ldots,x_{i}])[x_{i+1}]$ of the form

$$q_i(x_1,\ldots,x_{i+1})=\prod_{l=0}^{\mathrm{deg}_{x_{i+1}}p_i}\frac{\partial^l p_i}{\partial x_{i+1}^l}(x_1,\ldots,x_{i+1}),$$
where $p_i\in(\mathbb{R}[x_1,\ldots,x_{i}])[x_{i+1}]$ is monic in $x_{i+1},$  such that
$$\tilde{\pi}_{i+1}(W_i)=\mathcal{Z}(q_i).$$
For $d=i+1,\ldots,n-1,$
$$\hspace*{-42mm}(4_{i,d})\hspace*{41mm}  (\tilde{\pi}_d(W_i)\times\mathbb{R})\cap\tilde{\pi}_{d+1}(W_d)=\tilde{\pi}_{d+1}(W_i).$$
For $d=i+1,\ldots,n-1,$
$$\hspace*{-44mm}(5_{i,d})\hspace*{42mm}\tilde{\pi}_d(W_{i+1})\setminus\tilde{\pi}_d(W_{i})\subset\mathrm{Reg}(\tilde{\pi}_d(W_{i+1}))$$ and
$\tilde{\pi}_{d+1}(W_{i+1})$ is locally the union of pairwise disjoint graphs of analytic functions over
$\tilde{\pi}_d(W_{i+1})\setminus\tilde{\pi}_d(W_{i}).$
\end{proposition}

\begin{remark} \label{algebraicWi} The conditions $(3_i)$, $(4_{i,d}),$ $d=i+1,\ldots,n-1,$  of
Proposition~\ref{formulation},
in combination with induction on $d$ and Lemma \ref{localsur}, imply that for every
 $d=i+1,\ldots,n,$ the set $\tilde{\pi}_d(W_i)$ is an algebraic subset of $\mathbb{R}^d$ of pure dimension $i.$
\end{remark}

\begin{remark} \label{surjectglob} The condition $(3_i)$ of Proposition \ref{formulation} and Lemma \ref{localsur} imply
that  $$\tilde{\pi}_i(W_i)=\mathbb{R}^i.$$
\end{remark}
\textit{Proof of Remark \ref{surjectglob}.} We have $$\tilde{\pi}_{i}(W_i)=\pi_i(\mathcal{Z}(q_i)).$$ Since $p_i$
is a monic polynomial, the map $\pi_i|_{\mathcal{Z}(q_i)}$ is proper and therefore $\pi_i({\mathcal{Z}(q_i)})$ is
a closed subset of $\mathbb{R}^i.$ Now, Lemma \ref{localsur} immediately implies the assertion.\qed

\begin{remark} \label{globgraph} The condition $(5_{i,d})$ of Proposition \ref{formulation} implies that for every
connected component $B$ of $\tilde{\pi}_d(W_{i+1})\setminus\tilde{\pi}_d(W_i)$ there is a positive integer $k$
such that $(B\times\mathbb{R})\cap \tilde{\pi}_{d+1}(W_{i+1})$ is the union of $k$ pairwise disjoint graphs of
analytic functions $B\rightarrow\mathbb{R}.$ Moreover, by $(4_{i+1,d})$ we have
$$(B\times\mathbb{R})\cap \tilde{\pi}_{d+1}(W_{i+1})=(B\times\mathbb{R})\cap \tilde{\pi}_{d+1}(W_{d}).$$
\end{remark}

\begin{remark}\label{connectedcomps} The condition $(2_{i,j})$ of Proposition \ref{formulation} implies that either $(W_{i+1}\cap S_j)\setminus W_i=\emptyset$ or every connected component of
$(W_{i+1}\cap S_j)\setminus W_i$ is a connected component of $W_{i+1}\setminus W_i.$
\end{remark}
\textit{Proof of Remark \ref{connectedcomps}}. Let $A\neq\emptyset$ be a connected component of $(W_{i+1}\cap
S_j)\setminus W_i.$ Suppose that there is a connected set $B\subset W_{i+1}\setminus W_i$ such that $A$ is a
proper subset of $B.$ Then $B$ has common points both with $W_{i+1}\cap S_j$ and $W_{i+1}\setminus S_j.$ It
follows that
$$B\cap\partial_{W_{i+1}}(W_{i+1}\cap S_j)\neq\emptyset.$$ This contradicts $(2_{i,j})$ of Proposition \ref{formulation},
as $B\cap W_i=\emptyset.$\qed\\
\\
\textit{\textbf{Proof of Proposition \ref{formulation}.}} We construct $W_i$ using descending induction on $i.$
Let $i=n-1.$

Let $\tilde{W}_{n-1}$ be an algebraic hypersurface in $\mathbb{R}^n$ such that $
\partial_{\mathbb{R}^n}S_j\subset\tilde{W}_{n-1}$ for every $j\in J.$ After a linear change of coordinates we have
$$\tilde{W}_{n-1}=\mathcal{Z}(p_{n-1}),$$ where $p_{n-1}\in(\mathbb{R}[x_1,\ldots,x_{n-1}])[x_n]$ is monic in
$x_n.$ Let $r_{n-1}$ be the degree of $p_{n-1}$ in $x_n.$ Define
$$q_{n-1}(x_1,\ldots,x_n):=\prod_{l=0}^{r_{n-1}}\frac{\partial^l p_{n-1}}{\partial x_n^l}(x_1,\ldots,x_n) \mbox{
and }W_{n-1}:=\mathcal{Z}(q_{n-1}).$$ Clearly, $(1_{n-1}), (3_{n-1})$ and $(2_{n-1,j})$ of Proposition
\ref{formulation} hold true for every $j\in J$ with $W_{n-1}$ as above and $W_n=\mathbb{R}^n.$ The conditions
$(4_{n-1,d}), (5_{n-1,d})$ are trivially satisfied as there is no $d$ with $n\leq d\leq n-1.$

Let $i<n-1.$ Assume that  for $l=i+1,\ldots, n-1$ there exists an algebraic subset $W_{l}$ of $\mathbb{R}^n$
satisfying $(1_{l}), (2_{l,j}), (3_{l}),$ $(4_{l,d}),$ $(5_{l,d}),$ for $j\in J$ and $d=l+1,\ldots,n-1.$ We will
construct $W_i$ satisfying $(1_{i}), (2_{i,j}), (3_{i}),$ $(4_{i,d}),$ $(5_{i,d}),$ for $j\in J$ and
$d=i+1,\ldots,n-1.$\vspace*{3mm}

\noindent\textbf{Claim A. }\textit{For $d=i+1,\ldots, n,$ there is  an algebraic set
$D^d_{i}\subset\mathbb{R}^d\approx\mathbb{R}^{i+1}\times\mathbb{R}^{d-i-1}$ of dimension at most $i,$ contained in
$\tilde{\pi}_d(W_{i+1})$ such that:}\vspace*{2mm}

\noindent\textit{$(\sharp)$ $ \tilde{\pi}_{d+1}(W_{i+1})$ is locally the union of pairwise disjoint graphs of
analytic functions over $\tilde{\pi}_{d}(W_{i+1})\setminus D^d_{i},$ { for } $d=i+1,\ldots,n-1,$}\vspace*{2mm}

\noindent$(\dag)$ $\mathrm{Sing}(\tilde{\pi}_{d}(W_{i+1}))\subset D^d_{i},$\textit{ for }
$d=i+1,\ldots,n,$\vspace*{2mm}

 \noindent$(\flat)$
$\partial_{W_{i+1}}(S_j\cap W_{i+1})\subset D^n_{i},$\textit{ for }$j\in J.$\vspace*{4mm}

\noindent{\textit{Proof of Claim A.}} By Remark \ref{algebraicWi}, for $d=i+2,\ldots,n,$ the set
$\tilde{\pi}_d(W_{i+1})$ is an algebraic subset of $\mathbb{R}^d$ of pure dimension $i+1.$ By $(3_{l}),$ for
$l=i+1,\ldots,n-1,$ the set $\tilde{\pi}_{l+1}(W_{l})$ is an algebraic subset of $\mathbb{R}^{l+1}$ with proper
projection of finite fibers to $\mathbb{R}^l.$ Combining these two facts with $(4_{i+1,l})$ for $l=i+2,\ldots,n-1$
we conclude that $\tilde{\pi}_{d+1}(W_{i+1})$ is an algebraic subset of $\mathbb{R}^{d+1}$ of pure dimension
$i+1,$ with proper surjective projection of finite fibers to $\tilde{\pi}_{d}(W_{i+1})\subset\mathbb{R}^d,$ for
$d=i+1,\ldots,n-1.$ This implies existence of $D^d_i\subset\tilde{\pi}_d(W_{i+1})$ of dimension $i,$
$d=i+1,\ldots,n-1,$ satisfying $(\sharp).$

Since $\mathrm{Sing}(\tilde{\pi}_{d}(W_{i+1}))$ is an algebraic subset of $\mathbb{R}^{d}$ of dimension at most
$i,$ we may assume, enlarging $D^d_i,$ that $(\dag)$ is satisfied.

Finally, $\mathrm{dim}(W_{i+1})=i+1,$ so the dimension of $\partial_{W_{i+1}}(S_j\cap W_{i+1})$ is at most $i.$
Hence, enlarging $D^n_i$ we may assume that $(\flat)$ is also satisfied.\qed\vspace*{2mm}

We return to the proof of Proposition \ref{formulation}. Let $\tilde{D}^d_i$ be the image of the projection of
$D^d_i$ (provided by Claim A) to $\mathbb{R}^{i+1},$ for $d=i+1,\ldots,n.$ Define $\Theta_i$ to be the Zariski
closure in $\mathbb{R}^{i+1}$ of the union of all $\tilde{D}^d_i.$ Clearly, $\Theta_i$ is an algebraic subset of
$\mathbb{R}^{i+1}$ of dimension at most $i.$ Let $p_i\in(\mathbb{R}[x_1,\ldots,x_i])[x_{i+1}]$ such that
$\mathcal{Z}(p_i)=\Theta_i.$ We may assume, applying a linear change of coordinates in $\mathbb{R}^{i+1}$, that
$p_i$ is monic in $x_{i+1}.$

Define $$q_i(x_1,\ldots,x_{i+1}):=\prod_{j=0}^{\mathrm{deg}_{x_{i+1}}(p_i)}\frac{\partial^j p_i}{\partial
x_{i+1}^j}(x_1,\ldots,x_{i+1})\mbox{ and }C^{i+1}_i:=\mathcal{Z}(q_i)\subset\mathbb{R}^{i+1}.$$ Next, for
$d=i+1,\ldots,n-1,$ define inductively
$$C^{d+1}_i:=(C^d_i\times\mathbb{R})\cap\tilde{\pi}_{d+1}(W_d)\subset\mathbb{R}^{d+1}.$$
Finally, define $$W_i:=C^n_i.$$

It remains to show that $W_i$ enjoys the required properties.\vspace*{2mm}

\noindent\textbf{Claim B.} $$ C_i^{d}=\tilde{\pi}_d(W_i), \mbox{ for }d=i+1,\ldots,n.$$ \textit{Proof.} Descending
induction on $d.$ For $d=n,$ the equality follows by definition of $W_i.$

Assume that $i+1\leq d<n$ and that $ C_i^{d+1}=\tilde{\pi}_{d+1}(W_i).$ By definition of $C_i^{d+1}$ we have
$$\tilde{\pi}_{d+1}(W_i)=(C^d_i\times\mathbb{R})\cap\tilde{\pi}_{d+1}(W_d).$$
By $(3_d)$ and Remark \ref{surjectglob}, we have $\pi_d(\tilde{\pi}_{d+1}(W_d))=\mathbb{R}^d.$ Consequently,
$$\tilde{\pi}_{d}(W_i)=\pi_d(\tilde{\pi}_{d+1}(W_i))=\pi_d((C^d_i\times\mathbb{R})\cap\tilde{\pi}_{d+1}(W_d))=C^d_i.$$
The proof of Claim B is complete. \qed\vspace*{2mm}


By Claim B and the definition of $C^{i+1}_i$ we have $$\tilde{\pi}_{i+1}(W_i)=C^{i+1}_i=\mathcal{Z}(q_i),$$ hence
$W_i$ satisfies $(3_i).$

Next, Claim B and the inductive definition of $C^{d+1}_i$ immediately imply  $(4_{i,d})$ for $d=i+1,\ldots,n-1.$

The following fact will be useful in the proof of $(1_i), (2_{i,j}),$ $(5_{i,d}),$ for
$d=i+1,\ldots,n-~1.$\vspace*{2mm}

\noindent\textbf{Claim C.}
$$D^d_i\subset\tilde{\pi}_{d}(W_i), \mbox{ for }d=i+1,\ldots,n.$$
\textit{Proof.} Fix $d.$ For $l<d,$ define $\hat{\pi}_{l,d}:\mathbb{R}^d\rightarrow\mathbb{R}^l$ by
$\hat{\pi}_{l,d}:=\pi_l\circ\ldots\circ\pi_{d-1}.$ Put $\hat{\pi}_{d,d}=\mathrm{id}_{\mathbb{R}^d}$. By
construction and Claim B,

$$\hat{\pi}_{i+1,d}(D^d_i)\subset\Theta_i\subset
C^{i+1}_i=\tilde{\pi}_{i+1}(W_i)=\hat{\pi}_{i+1,d}( \tilde{\pi}_{d}(W_i)).$$  Now assume that
$$\hat{\pi}_{l,d}(D^d_i)\subset\hat{\pi}_{l,d}( \tilde{\pi}_{d}(W_i))=\tilde{\pi}_l(W_i)$$ for some $l\in\{i+1,\ldots,d-1\}.$ Then, by
$(4_{i,l}),$ we get
$$\hspace*{-43mm}(\ddag)\hspace*{43mm}(\hat{\pi}_{l,d}(D^d_i)\times\mathbb{R})\cap\tilde{\pi}_{l+1}(W_l)\subset\tilde{\pi}_{l+1}(W_i).$$
By Claim A and $(1_{i+1}),\ldots,(1_{l-1})$, we have

 $$D_i^d\subset \tilde{\pi}_d(W_{i+1})\subset\tilde{\pi}_d(W_l),$$ so

$$\hat{\pi}_{l+1,d}(D_i^d)\subset\hat{\pi}_{l+1,d}(\tilde{\pi}_d(W_l))=\tilde{\pi}_{l+1}(W_l),$$
which in conjunction with $(\ddag)$ gives
$$\hat{\pi}_{l+1,d}(D_i^d)\subset\tilde{\pi}_{l+1}(W_i)=\hat{\pi}_{l+1,d}(\tilde{\pi}_d(W_i)).$$
Using induction on $l,$ we get
$$D^d_i=\hat{\pi}_{d,d}(D_i^d)\subset\hat{\pi}_{d,d}(\tilde{\pi}_d(W_i))=\tilde{\pi}_d(W_i).$$
The proof of Claim C is complete.\qed\vspace*{2mm}

Let us check $(1_i).$ Using induction on $d,$ we check that $\tilde{\pi}_d(W_i)\subset \tilde{\pi}_d(W_{i+1}),$
for $d=i+1,\ldots,n.$ For $d=i+1$ this is obvious as $\tilde{\pi}_{i+1}(W_{i+1})=\mathbb{R}^{i+1}.$ The condition
$(4_{i,i+1})$ gives $\tilde{\pi}_{i+2}(W_i)\subset \tilde{\pi}_{i+2}(W_{i+1}).$ The fact that
$\tilde{\pi}_{d+1}(W_i)\subset \tilde{\pi}_{d+1}(W_{i+1})$ immediately follows from $\tilde{\pi}_d(W_i)\subset
\tilde{\pi}_d(W_{i+1})$ by $(4_{i,d}), (4_{i+1,d}),$ for $d\geq i+2,$ finishing the inductive procedure. In
particular, for $d=n,$ we have $W_i\subset W_{i+1}.$

By $(4_{i,d})$ and $(3_d),$ for $d=i+1,\ldots,n-1,$ we have
$$\mathrm{dim}(\tilde{\pi}_d(W_i))=\mathrm{dim}(\tilde{\pi}_{d+1}(W_i)).$$ Hence, by $(3_i)$
$$ i=\mathrm{dim}(\tilde{\pi}_{i+1}(W_i))=\mathrm{dim}(\tilde{\pi}_{n}(W_i))=\mathrm{dim}(W_i).$$

Finally, to see that $\mathrm{Sing}(W_{i+1})\subset W_i,$  it is sufficient to observe, in view of $(\dag),$ that
$D^n_i\subset W_i,$ which follows by Claim C and proves $(1_i).$

The fact that $D_i^n\subset W_i$ also proves $(2_{i,j}),$ for $j\in J,$ in view of $(\flat).$

It remains to check that $(5_{i,d})$ holds true for $d=i+1,\ldots,n-1.$ The inclusion in $(5_{i,d})$ follows from
$\mathrm{Sing}(\tilde{\pi}_d(W_{i+1}))\subset \tilde{\pi}_d(W_{i})$ which is a consequence of Claim~C in
conjunction with $(\dag).$ The remaining part of $(5_{i,d})$ follows by Claim~C  and $(\sharp).$ The proof of
Proposition \ref{formulation} is complete.\qed
\subsection{Cell decomposition}

The notion of a cell decomposition is a classical topic in the theory of o-minimal structures (cf. \cite{VdD}).
Cells which appear in the present paper are cells in the sense of o-minimality, but satisfy additional restrictive
conditions which allow us to control the behavior of piecewise-regular functions.

We adopt the following notation. Let $E\subset\mathbb{R}^{d-1}$ and let $f_1:E\rightarrow\mathbb{R}$ or $f_1$
constantly $-\infty$ on $E,$ and $f_2:E\rightarrow\mathbb{R}$ or $f_2$ constantly $+\infty$ on $E.$  Put
$$(f_1,f_2)_E:=\{(x,y)\in E\times\mathbb{R}: f_1(x)<y< f_2(x)\},$$
$$[f_1,f_2]_E:=\{(x,y)\in E\times\mathbb{R}: f_1(x)\leq y\leq f_2(x)\}.$$
In particular, $[f_1,f_1]_E=\mathrm{graph}(f_1)$ if $f_1$ is real valued, and $[f_1,f_1]_E=\emptyset$ if $f_1$ is
$+/-\infty.$ We write $f_1<f_2$ to express that $f_1(x)<f_2(x)$ for every $x\in E.$

For $i=0,\ldots,n-1,$ let $\pi_i:\mathbb{R}^i\times\mathbb{R}\rightarrow\mathbb{R}^i$ be the natural projection.
For $d=1,\ldots,n-1,$ define $\tilde{\pi}_d:\mathbb{R}^n\rightarrow\mathbb{R}^d,$
$\tilde{\pi}_d:=\pi_d\circ\pi_{d-1}\circ\cdots\circ\pi_{n-1}.$ Put $\tilde{\pi}_n:=\mathrm{Id}_{\mathbb{R}^n}.$
For $p\in\mathbb{R}[x_1,\ldots,x_{i+1}],$ let $\mathrm{deg}_{x_{i+1}}p$ be the degree of $p$ in $x_{i+1}.$

\subsubsection{Construction of cells}\label{constructofcells}

Let $W_0,\ldots, W_{n}$ be real algebraic subsets  of $\mathbb{R}^n,$ $W_n=\mathbb{R}^n,$ such that {for} every
$i=0,\ldots,n-1$ and $d=i+1,\ldots,n-1$ satisfied are: $(1_i),$ $(3_{i}),$ $(4_{i,d}),$ $(5_{i,d})$ (contained in
Proposition \ref{formulation}). Put $W_{-1}=\emptyset.$

For every $d=1,\ldots,n$ we will construct a finite family $\mathcal{K}_d$ of cells in $\mathbb{R}^d$ whose union
gives $\mathbb{R}^d$ such that each cell is contained in $\tilde{\pi}_d(W_{i+1})\setminus\tilde{\pi}_d(W_i)$ for
some $i\in\{-1,0,\ldots,d-1\}.$ The construction is carried out by induction on $d.$

For $d=1,$ each point in $\tilde{\pi}_1(W_0)$ is a by definition a cell of dimension $0$. Moreover, each connected
component of $\tilde{\pi}_1(W_1)\setminus\tilde{\pi}_1(W_0)=\mathbb{R}\setminus\tilde{\pi}_1(W_0)$ is by
definition a cell of dimension $1$.

Let $d>1.$ A subset $C\subset\mathbb{R}^{d-1}\times\mathbb{R}$ is a cell if there exist $i\in\{-1,0,\ldots,d-2\},$
a cell $B\subset\tilde{\pi}_{d-1}(W_{i+1})\setminus\tilde{\pi}_{d-1}(W_i)\subset\mathbb{R}^{d-1}$ and
either:\vspace*{2mm}\\
(a) functions $f<g$ on $B$ such that $f:B\rightarrow\mathbb{R}$ analytic or $f=-\infty,$ and
$g:B\rightarrow\mathbb{R}$ analytic or $g=+\infty$ with
$$[f,f]_B,[g,g]_B\subset\tilde{\pi}_d(W_{i+1}),$$ $$(f,g)_B\subset\mathbb{R}^d\setminus\tilde{\pi}_d(W_{{i+1}}),$$
and $C=(f,g)_B,$\vspace*{2mm}\\
or:\vspace*{2mm}

\noindent(b) an analytic function $f:B\rightarrow\mathbb{R}$ with
$$[f,f]_B\subset\tilde{\pi}_d(W_{i+1}),$$
and $C=[f,f]_B.$\vspace*{2mm}

In the case (a) put $\mathrm{dim}(C)=\mathrm{dim}(B)+1.$ In the case (b) put $\mathrm{dim}(C)=\mathrm{dim}(B).$

\subsubsection{Properties of cells}

\begin{remark}\label{Obscell} If $C=(f,g)_B$ is a cell in $\mathbb{R}^d,$ then by $(4_{i+1,d}),$ we have
$C \subset\mathbb{R}^d\setminus\tilde{\pi}_d(W_{{d-1}}).$ Moreover,
$$[f,f]_B,[g,g]_B\subset\tilde{\pi}_d(W_{d-1}),$$ as $W_{i+1}\subset W_{d-1}.$ Thus $C$ is a region between two
consecutive graphs  over $B$ which are contained in $\tilde{\pi}_d(W_{d-1}).$

\end{remark}

Let us list basic properties of cells.

\begin{proposition}\label{basicprocell}~\vspace*{2mm}

\noindent(a) For each cell  $C$ in $\mathbb{R}^d$ there is $i\in \{-1,0,\ldots,d-1\}$ such that
$$C\subset\tilde{\pi}_d(W_{i+1})\setminus\tilde{\pi}_d(W_i).$$

In particular, $(C\times\mathbb{R})\cap\tilde{\pi}_{d+1}(W_d)$ is the union of pairwise disjoint graphs of
analytic

maps $C\rightarrow\mathbb{R}.$

\noindent(b) Every cell is connected.

\noindent(c) The union of all cells in $\mathbb{R}^d$ equals $\mathbb{R}^d.$\\
(d) If $A, B$ are cells in $\mathbb{R}^{d}$ with $A\neq B$, then $A\cap B=\emptyset.$\\
(e) If $C$ is a cell in $\mathbb{R}^d$ of dimension $d,$ then $C$ is an open subset of $\mathbb{R}^d.$\\
(f) Let $\{S_j\}_{j\in J}$ be a finite family of semialgebraic subsets of $\mathbb{R}^n$ satisfying $(2_{i,j})$
for all

$i,j.$ Then for every $j$ and every cell $C$ in $\mathbb{R}^n$ either $C\subset S_j$ or $C\cap S_j=\emptyset.$

\end{proposition}
\textit{Proof.} (a) Induction on $d.$ For $d=1$ there is nothing to prove. If $d>1,$ then either $C=(f,g)_B$ or
$C=[f,f]_B.$ In the first case by Remark \ref{Obscell}, we have $$C
\subset\mathbb{R}^d\setminus\tilde{\pi}_d(W_{{d-1}})=\tilde{\pi}_{d}(W_d)\setminus\tilde{\pi}_d(W_{{d-1}}).$$

In the second case, by definition, $C=[f,f]_B\subset\tilde{\pi}_d(W_{i+1})$ and
$B\cap\tilde{\pi}_{d-1}(W_i)=\emptyset.$ Therefore we have $C\subset
\tilde{\pi}_d(W_{i+1})\setminus\tilde{\pi}_d(W_i).$ The second part by Remark \ref{globgraph} and (b).

(b) and (e). Immediate by definition of cell and induction on $d.$

(c) and (d). Immediate by definition of cell, induction on $d,$ (a) and Remark \ref{globgraph}.

(f). Let $C$ be a cell in $\mathbb{R}^n.$ By (a), there is $i$ such that $C\subset W_{i+1}\setminus W_i.$ Since
$C$ is connected it is contained in a connected component of $W_{i+1}\setminus W_i.$ If $C\cap S_j\neq\emptyset,$
then there is a connected component $E$ of $(W_{i+1}\cap S_j)\setminus W_i$ such that $E\cap C\neq\emptyset.$ By
Remark \ref{connectedcomps}, $E$ is a connected component of $W_{i+1}\setminus W_i,$ hence $C\subset E.$
Consequently, $C\subset S_j$ as $E\subset S_j.$\qed\vspace*{2mm}

\noindent\textbf{Definition.} \textit{A closed cell in $\mathbb{R}^d$ of dimension $l$ is the closure of a cell of
dimension $l$ in the Euclidean topology of $\mathbb{R}^d.$}\vspace*{2mm}

\noindent\textbf{Remark.} \textit{Let $A$ be a cell in $\mathbb{R}^d.$ Induction on $d$ shows that
$\overline{A}\setminus A$ is a closed subset of $\mathbb{R}^d.$}\vspace*{2mm}

Clearly, every closed cell in $\mathbb{R}$ is either a point or $[a,b]$ or $[a,\infty)$ or $(-\infty,b]$ for some
$a<b.$

For any $x,y\in\mathbb{R}^d$ let $\overline{xy}$ denote the closed segment connecting $x$ and $y.$

\begin{lemma}\label{closlem}Let $d>i$ and $E\subset\tilde{\pi}_{d}(W_{i+1})\setminus\tilde{\pi}_{d}(W_{i})$ be a cell. Then
there are a homeomorphism
$\Theta:\overline{E}\rightarrow\Theta(\overline{E})=\overline{\Theta(E)}\subset\mathbb{R}^d$ and $x_E\in\Theta(E)$
such that:\vspace*{2mm}

\noindent(a) for each $y\in\overline{\Theta(E)},$ we have $\overline{x_Ey}\subset\overline{\Theta(E)}$ and $\overline{x_Ey}\cap(\overline{\Theta(E)}\setminus\Theta(E))\subset\{y\}$,\\
(b) if $y\in\overline{\Theta(E)},$ $y\neq x_E,$ and $\overline{x_Ey}$ is a maximal segment starting at $x_E$

with $\overline{x_Ey}\subset\overline{\Theta(E)}$, then $y\in\overline{\Theta(E)}\setminus\Theta(E).$\\
(c) each continuous function $f:E\rightarrow\mathbb{R}$ with $[f,f]_E\subset\tilde{\pi}_{d+1}(W_{i+1})$ extends
continuously to $\overline{E}\rightarrow\mathbb{R}.$
\end{lemma}
\textit{Proof.} Induction on $d.$ Let $d=1.$ If $E$ is a point, there is nothing to prove. If $E$ is an interval
then take $\Theta=\mathrm{id}_{\mathbb{R}}.$ Then $(a), (b)$ are trivial. As for $(c),$ we have
$$[f,f]_E\subset\tilde{\pi}_{d+1}(W_d)=\mathcal{Z}(q_d).$$ Since the leading coefficient of $q_d$ is a
nonzero constant, by the continuity of roots of $q_d$ the function $f$ extends continuously to
$\overline{E}\rightarrow\mathbb{R}.$

Let $d>1.$ We have either $E=(g,h)_F$ of $E=[g,g]_F$ for some cell $F$ in $\mathbb{R}^{d-1}$ and $g<h$ (possibly
$g=-\infty$ or $h=+\infty$ in the first case).

By induction hypothesis, there are a homeomorphism
$\Xi:\overline{F}\rightarrow\Xi(\overline{F})=\overline{\Xi(F)}\subset\mathbb{R}^{d-1}$ and $x_F\in\Xi(F)$ such
that:\vspace*{2mm}

\noindent(u) for each $y\in\overline{\Xi(F)},$ we have $\overline{x_Fy}\subset\overline{\Xi(F)}$
and $\overline{x_Fy}\cap(\overline{\Xi(F)}\setminus\Xi(F))\subset\{y\}$,\\
(v) if $y\in\overline{\Xi(F)},$ $y\neq x_F,$ and $\overline{x_Fy}$ is a maximal segment starting at $x_F$
with $\overline{x_Fy}\subset\overline{\Xi(F)}$, then $y\in\overline{\Xi(F)}\setminus\Xi(F).$\\

First consider the case where $E=(g,h)_F.$ By induction hypothesis the functions $g,h$ have continuous extensions
defined on $\overline{F}$ also denoted by $g, h.$ Then we have,
$$\overline{E}=[{g},{h}]_{\overline{F}}\mbox{ and }(\Xi\times\mathrm{id}_{\mathbb{R}})({E})=(\tilde{g},\tilde{h})_{\Xi({F})}\mbox{ and
}(\Xi\times\mathrm{id}_{\mathbb{R}})(\overline{E}) =[\tilde{g},\tilde{h}]_{\Xi(\overline{F})},$$ where
$\tilde{g}={g}\circ\Xi^{-1}|_{\Xi(\overline{F})}, \tilde{h}={h}\circ\Xi^{-1}|_{\Xi(\overline{F})}.$

If both $\tilde{g}=-\infty$ and $\tilde{h}=+\infty,$ we define $\Theta=(\Xi\times\mathrm{id}_{\mathbb{R}}),$ which
clearly satisfies (a) and (b). If $\tilde{g}=-\infty$ and $\tilde{h}<+\infty,$ then the map
$$\gamma:[\tilde{g},\tilde{h}]_{\Xi(\overline{F})}\rightarrow (-\infty,0]_{\Xi(\overline{F})}$$ given by the formula
$\gamma(x,t)=(x,t-\tilde{h}(x))$ is a homeomorphism. In this case define $\Theta=\gamma\circ
(\Xi\times\mathrm{id}_{\mathbb{R}})|_{\overline{E}},$ which clearly satisfies (a) and (b). The case
$\tilde{g}>-\infty$ and $\tilde{h}=+\infty$ is very similar.

Now assume that $\tilde{g}, \tilde{h}:\Xi(\overline{F})\rightarrow\mathbb{R}.$ Elementary constructions show that
there is a homeomorphism
$$\alpha:[\tilde{g},\tilde{h}]_{\Xi(\overline{F})}\rightarrow[0,\mathrm{min}\{1,\tilde{h}-\tilde{g}\}]_{\Xi(\overline{F})}$$
with $\alpha((\tilde{g},\tilde{h})_{\Xi({F})})=(0,\mathrm{min}\{1,\tilde{h}-\tilde{g}\})_{\Xi({F})}.$ Indeed, for
every fixed $x\in\Xi(\overline{F})$ it is sufficient to define $\alpha(x,\cdot)$ by an affine-linear interpolation
formula transforming $\tilde{g}(x), \tilde{h}(x)$ to $0, \mathrm{min}\{1, \tilde{h}(x)-\tilde{g}(x)\},$
respectively. If $\tilde{h}(x)=\tilde{g}(x),$ then $\tilde{h}(x)$ is mapped to $0.$

Define a function $u:\Xi(\overline{F})\rightarrow [0,1]$ as follows. Put $u(x_F)=1.$ For any  ray
$T\subset\Xi(\overline{F})$ maximal with respect to the inclusion, starting from $x_F$, proceed as follows. If $T$
is unbounded, put $u|_T=1.$ If $T$ is bounded, then $T=\overline{x_F\hat{x}},$ where $\hat{x}$ is the unique
intersection point of $T$ and $\overline{\Xi(F)}\setminus\Xi(F).$ Then for every $x\in T$, there is $t_x\in [0,1]$
with $x=x_F+t_x(\hat{x}-x_F).$ Define
$$u(x):=1-t_x+t_x\mathrm{min}\{1,\tilde{h}(\hat{x})-\tilde{g}(\hat{x})\}.$$

Since $\Xi(\overline{F})$ is covered by rays starting from $x_F,$ $u$ is a well defined function on
$\Xi(\overline{F}).$ It is easy to see that $u$ is continuous and that $u(x)=0$ if and only if
$\mathrm{min}\{1,\tilde{h}(x)-\tilde{g}(x)\}=0.$ Therefore the map

$$\beta:[0,\mathrm{min}\{1,\tilde{h}-\tilde{g}\}]_{\Xi(\overline{F})}\rightarrow[0,u]_{\Xi(\overline{F})}$$
defined for every $x\in\Xi(\overline{F})$ by affine-linear interpolation formula mapping $0,
\mathrm{min}\{1,\tilde{h}(x)-\tilde{g}(x)\}$ to $0, u(x),$ respectively, is a homeomorphism. Then for
$\Theta=\beta\circ\alpha\circ(\Xi\times\mathrm{id}_{\mathbb{R}})|_{\overline{E}},$ we have
$\Theta(E)=(0,u)_{\Xi(F)}$ and $\Theta(\overline{E})=\overline{\Theta(E)}=[0,u]_{\Xi(\overline{F})}.$

Note that for every ray $T\subset\Xi(\overline{F})$ starting from $x_F$ the graph of the restriction $u|_{T}$ is a
segment. It is easy to see that (a), (b) are satisfied with $x_E=(x_F,\frac{1}{2}).$

Now assume that $E=[g,g]_F.$ By induction hypothesis the function $g$ has a continuous extension
$\overline{F}\rightarrow\mathbb{R},$ denoted by $g.$ Then we have,
$$\overline{E}=[{g},{g}]_{\overline{F}}\mbox{ and }(\Xi\times\mathrm{id}_{\mathbb{R}})({E})
=[\tilde{g},\tilde{g}]_{\Xi({F})}\mbox{ and }(\Xi\times\mathrm{id}_{\mathbb{R}})(\overline{E})
=[\tilde{g},\tilde{g}]_{\Xi(\overline{F})},$$ where $\tilde{g}={g}\circ\Xi^{-1}|_{\Xi(\overline{F})},
\tilde{h}={h}\circ\Xi^{-1}|_{\Xi(\overline{F})}.$ Clearly, $[\tilde{g},\tilde{g}]_{\Xi(\overline{F})}$ is
homeomorphic to $[0,0]_{\Xi(\overline{F})}$ which satisfies both (a) and (b) with $x_E=(x_F,0)$. Thus we have
proved that there is $\Theta$ satisfying (a), (b).

Now let $f:E\rightarrow\mathbb{R}$ be a continuous function as in (c). Then
$f\circ(\Theta|_E)^{-1}:\Theta(E)\rightarrow\mathbb{R}$ is a continuous function. Since
$$[f,f]_E\subset\tilde{\pi}_{d+1}(W_d)=\mathcal{Z}(q_d),$$ we have

$$[f\circ(\Theta|_E)^{-1},f\circ(\Theta|_E)^{-1}]_{\Theta(E)}\subset\mathcal{Z}(q_d(\Theta^{-1},x_{d+1})).$$

Fix any $c\in\overline{\Theta(E)}\setminus\Theta(E).$ By (a), we have $\overline{x_Ec}\subset\overline{\Theta(E)}$
and $\overline{x_Ec}\setminus\{c\}\subset \Theta(E).$  Continuity of roots of $q_d(\Theta^{-1},x_{d+1})$ implies
that  there is the unique value which assigned to $c$ extends the function
$f\circ(\Theta^{-1}|_{(\overline{x_Ec}\setminus\{c\})})$ continuously to $\overline{x_Ec}\rightarrow\mathbb{R}.$

Assigning value to each $c\in\overline{\Theta(E)}\setminus\Theta(E)$ as in the previous paragraph gives the
continuous extension $g:\overline{\Theta(E)}\rightarrow\mathbb{R}$ of
$f\circ(\Theta|_E)^{-1}:\Theta(E)\rightarrow\mathbb{R}$ (again by continuity of roots of
$q_d(\Theta^{-1},x_{d+1})$). Since $\overline{\Theta(E)}=\Theta(\overline{E})$ the function
$g\circ\Theta:\overline{E}\rightarrow\mathbb{R}$ extends continuously $f:E\rightarrow\mathbb{R}.$\qed

\begin{lemma}\label{projelem}Let $C$ be a cell in $\mathbb{R}^d\times\mathbb{R}.$
Let $\pi:\mathbb{R}^d\times\mathbb{R}\rightarrow\mathbb{R}^d$ be the natural projection. Then
$$\overline{\pi(C)}=\pi(\overline{C}).$$
\end{lemma}
\textit{Proof.} If $C=[f,f]_D,$ where $D$ is a cell in $\mathbb{R}^d,$ then, by Lemma~\ref{closlem}(c), $f$ has a
continuous extension $\bar{f}:\overline{D}\rightarrow\mathbb{R}.$ Then
$\overline{C}=[\bar{f},\bar{f}]_{\overline{D}}$ and the assertion follows immediately.

Let $C=(f,g)_D,$ where $D$ is a cell in $\mathbb{R}^d.$ By Lemma~\ref{closlem}(c), $f,g$ have continuous
extensions $\bar{f},\bar{g},$ respectively, defined on $\overline{D}$. Then
$\overline{C}=[\bar{f},\bar{g}]_{\overline{D}}$ and the assertion follows immediately.\qed

\begin{lemma}\label{closedcell} Each closed cell in $\mathbb{R}^d=\mathbb{R}^{d-1}\times\mathbb{R}$ is of one of the following
forms:\vspace*{2mm}

\noindent(a) $[f,f]_D\subset\tilde{\pi}_{d}(W_{d-1}),$ where $D=\overline{C},$ and $C$ is a cell in
$\mathbb{R}^{d-1},$ and $f:D\rightarrow\mathbb{R}$ is a continuous function.\vspace*{2mm}

\noindent(b) $[f,g]_D$ with $[f,f]_D, [g,g]_D\subset \tilde{\pi}_{d}(W_{d-1}),$ where $D=\overline{C}$ and $C$ is
a cell in $\mathbb{R}^{d-1}.$ Moreover, $f:D\rightarrow\mathbb{R}$ is a continuous function or $f=-\infty,$ and
$g:D\rightarrow\mathbb{R}$ is a continuous function or $g=+\infty,$ and $f|_C<g|_C,$ and
$$(f|_C,g|_C)_C\subset\mathbb{R}^{d}\setminus\tilde{\pi}_d(W_{d-1}).$$

\noindent Furthermore, in case (b) either $f>-\infty$ or $g<+\infty.$
\end{lemma}
\textit{Proof.} Let $E$ be a closed cell in $\mathbb{R}^{d-1}\times\mathbb{R}.$ Then $E=\overline{F},$ where $F$
is a cell, that is, $F=[u,u]_C$ or $F=(u,v)_C,$ where $C$ is a cell in $\mathbb{R}^{d-1}.$ By
Lemma~\ref{closlem}(c), there exist continuous extensions $f,g$ of $u, v,$ respectively, defined on
$D=\overline{C}.$ Then either $E=[f,f]_D$ or $E=[f,g]_D$ and (a), (b) follow by Section \ref{constructofcells}.

If in case (b) we had both $f=-\infty$ and $g=+\infty,$ then $(f,g)_C\cap\tilde{\pi}_d(W_{d-1})=\emptyset.$ This
contradicts Remark \ref{surjectglob}. \qed\vspace*{2mm}

\noindent\textbf{Definition. } Let $E$ be a closed cell in $\mathbb{R}.$ If $E$ is an interval, then
$A\subset\mathbb{R}$ is a \textit{face} of $E$ if $A=E$ or $A$ is one of the boundary points of $E.$ If $E$ is a
point, then $E$ is the only \textit{face} of $E.$

Let $d>1$ and $E$ a closed cell in $\mathbb{R}^d.$ Then $E=[f,f]_D$ as in Lemma \ref{closedcell}(a) or $E=[f,g]_D$
as in Lemma \ref{closedcell}(b), where $D$ is a closed cell in $\mathbb{R}^{d-1}.$ If $E=[f,f]_D,$ then $A$ is a
\textit{face} of $E$ is there is a face $G$ of $D$ such that $A=[f|_G,f|_G]_G.$ If $E=[f,g]_D,$ then $A$ is a
\textit{face} of $E$ if there is a face $G$ of $D$ such that $A=[f|_G,g|_G]_G$ or ($f>-\infty$ and
$A=[f|_G,f|_G]_G$) or ($g<+\infty$ and  $A=[g|_G,g|_G]_G$).

\begin{remark}\label{transitface} Let $A, B, C$ be closed cells in $\mathbb{R}^d.$ Then:\vspace*{2mm}

\noindent(a) if $A$ is a face of $B$ and $B$ is a face of $C,$ then $A$ is a face of $C,$\\
(b) if $A$ is a face of $B,$ then $A\subset{B}.$
\end{remark}

\begin{lemma}\label{facclosedcell} Let $G$ be  a closed cell in $\mathbb{R}^d.$ Then every face of $G$ is a closed
cell in $\mathbb{R}^d.$
\end{lemma}
\textit{Proof.} Induction on $d.$ If $d=1,$ there is nothing to prove.

Let $d>1.$ Assume that $G=[f,f]_D$ as in Lemma \ref{closedcell}(a), where $D$ is a closed cell in
$\mathbb{R}^{d-1}.$ Let $A$ be a face. Then there is a face $E$ of $D$ such that $A=[f|_E,f|_E]_E.$ By induction
hypothesis, $E$ is a closed cell. Then $A$ is also a closed cell.

Assume that $G=[f,g]_D,$ $f\neq g$ as in Lemma \ref{closedcell}(b), where $D$ is a closed cell in
$\mathbb{R}^{d-1}.$ Let $A$ be a face. If there is a face $E$ of $D$ with $A=[f|_E,f|_E]_E$ or $A=[g|_E,g|_E]_E$,
then $E$ is a closed cell by induction hypothesis and then $A$ is a closed cell.

If there is a face $E$ of $D$ such that $A=[f|_E,g|_E]_E$ then $E$ is a closed cell by induction hypothesis. That
is $E=\overline{F},$ where $F$ is a cell. In particular, $F\subset
\tilde{\pi}_{d-1}(W_i)\setminus\tilde{\pi}_{d-1}(W_{i-1}),$ for some $i.$ So either $f|_F<g|_F$ or $f|_F=g|_F$ by
Proposition \ref{basicprocell}(a). In the second case, $f|_E=g|_E,$ so $[f|_E,g|_E]_E$ is a closed cell.

In the first case, $(f|_F,g|_F)_F\cap\tilde{\pi}_d(W_{d-1})=\emptyset$ because otherwise, by Lemma \ref{localsur},
we had $(f|_C, g|_C)_C\cap\tilde{\pi}_d(W_{d-1})\neq\emptyset,$ where $C$ is a cell with $\overline{C}=D$
(contradicting Lemma \ref{closedcell}(b)). Therefore $[f|_E,g|_E]_E$ is a closed cell.\qed

\begin{lemma}\label{frontier}Let $C\subset\mathbb{R}^d$ be a cell of dimension $l$. Then $\overline{C}\setminus C$ is the union of
all faces of $\overline{C}$ of dimension $l-1.$ In particular, every closed cell is the union of cells.
\end{lemma}
\textit{Proof.} Induction on $d.$ If $d=1,$ there is nothing to prove.

Let $d>1.$ Let $C=[f,f]_D,$ where $D$ is a cell in $\mathbb{R}^{d-1}$ of dimension $l.$ By Lemma~\ref{closlem}(c),
$f$ has a continuous extension $\overline{D}\rightarrow\mathbb{R}$ also denoted by $f.$ Then
$\overline{C}=[f,f]_{\overline{D}}.$ After applying the induction hypothesis to $\overline{D}\setminus D$ we
obtain the assertion immediately.

Let $C=(f,g)_D,$ $f<g,$ where $D$ is a cell in $\mathbb{R}^{d-1}$ of dimension $l-1.$ By Lemma~\ref{closlem}(c),
$f, g $ have continuous extensions on $\overline{D}$ also denoted by $f, g.$ Then
$\overline{C}=[f,g]_{\overline{D}}.$ Let $\Lambda$ be the family of all faces of $\overline{D}$ of dimension
$l-2.$ By induction hypothesis, $\overline{D}\setminus D=\bigcup\Lambda.$ Let $\Lambda'\subset\Lambda$ be the
family of all faces $L$ such that $f|_L\neq g|_L.$ For $L\in\Lambda',$ the dimension of $[f|_L,g|_L]_L$ equals
$l-1.$ We have
$$\overline{C}\setminus C=[f,f]_{\overline{D}}\cup [g,g]_{\overline{D}}\cup\bigcup_{L\in\Lambda}[f|_L,g|_L]_L=[f,f]_{\overline{D}}\cup [g,g]_{\overline{D}}\cup\bigcup_{L\in\Lambda'}[f|_L,g|_L]_L.$$

The last statement follows from Lemma \ref{facclosedcell} and the induction on $l.$  The proof is complete.\qed
\begin{lemma}\label{intersection} Let $A, B\subset\mathbb{R}^d$ be closed cells. Then $A\cap B$ is the union
of closed cells.
\end{lemma}
\textit{Proof.} By Lemma \ref{frontier} and Proposition \ref{basicprocell}(d), $A\cap B$ is the union of cells
closures of which are contained in $A\cap B$ as $A$  and $B$ are closed.\qed

\subsection{Isolating families}

Let $\mathcal{K}_d,$ for $d=1,\ldots,n,$ be a (finite) family of cells in $\mathbb{R}^d$ as defined in
Section~\ref{constructofcells}. By a cell in $\mathbb{R}^d$ we mean an element of $\mathcal{K}_d.$ Recall that a
closed cell in $\mathbb{R}^d$ is the closure of a cell in the Euclidean topology of $\mathbb{R}^d.$

For every cell $E$ in $\mathbb{R}^d,$ put
$$\mathcal{I}(E)=\{C\subset\mathbb{R}^d: C \mbox{ is a cell of dimension } d \mbox{ and }E\subset\overline{C}\}.$$

The finite family $\mathcal{I}(E)$ will be called the \textit{isolating family} for $E.$ It can be equivalently
defined by induction on $d$ as presented below. More precisely, we give an inductive construction of the family
denoted by $\mathcal{J}(E)$ and discuss its correctness and properties and we show that
$\mathcal{J}(E)=\mathcal{I}(E)$.

Let $E$ be a cell in $\mathbb{R}.$ If $E=\{a\},$ then there are the uniquely determined $b,c\in\mathbb{R}$ or
$b=-\infty$ or $c=+\infty$ such that $(b,a),$ $(a,c)$ are cells in $\mathbb{R}.$ Let
$\mathcal{J}(E)=\{(b,a),(a,c)\}.$ If $E$ is an open interval in $\mathbb{R},$ then put $\mathcal{J}(E)=\{E\}.$

Let $d>1$ and $E$ be a cell in $\mathbb{R}^d.$ Then, by Section~\ref{constructofcells}, either $E=[f,f]_B$ or
$E=(f,g)_B$ for some cell $B$ in $\mathbb{R}^{d-1}.$\vspace*{2mm}

\noindent\textbf{Case 1.} Let $E=[f,f]_B,$ where $B$ is a cell in $\mathbb{R}^{d-1}$ and
$f:B\rightarrow\mathbb{R}$ is a continuous function with $[f,f]_B\subset\tilde{\pi}_d(W_{d-1}).$ Assume we already
have constructed $\mathcal{J}(B)=\{C_1,\ldots,C_s\}.$

Let $u, v$ such that either $u=-\infty$ or $u:B\rightarrow\mathbb{R}$ is continuous with
$[u,u]_B\subset\tilde{\pi}_d(W_{d-1}),$ and either $v=+\infty$ or $v:B\rightarrow\mathbb{R}$ is continuous with
$[v,v]_B\subset\tilde{\pi}_d(W_{d-1}).$ Moreover, assume that $u<f<v$ and
$$(u,f)_B\cap\tilde{\pi}_d(W_{d-1})=\emptyset=
(f,v)_B\cap\tilde{\pi}_d(W_{d-1}).$$

By Lemma \ref{closlem}(c), $u,v$ have continuous extensions on $\overline{B}$ also denoted by $u,v.$

For each $i\in\{1,\ldots,s\},$ let $N_i^1,\ldots,N_i^{j(i)}:\overline{C_i}\rightarrow\mathbb{R}$ be all continuous
maps with graphs contained in $\tilde{\pi}_d(W_{d-1}),$ whose restrictions to $B$ equal $v.$ If there are no such
maps, put $N_i=+\infty,$ otherwise, define $N_i:=\min\{N_i^1,\ldots, N_i^{j(i)}\}.$ Let
$M_i^1,\ldots,M_i^{k(i)}:\overline{C_i}\rightarrow\mathbb{R}$ be all continuous maps with graphs contained in
$\tilde{\pi}_d(W_{d-1}),$ whose restrictions to $B$ equal $u.$ If there are no such maps, put $M_i=-\infty,$
otherwise, define $M_i:=\max\{M_i^1,\ldots, M_i^{k(i)}\}.$ Let $\mathcal{J}(E,i)$ be the family of all connected
components of $(M_i,N_i)_{C_i}\setminus\tilde{\pi}_d(W_{d-1}).$

Define $\mathcal{J}(E)=\bigcup_{i=1}^s\mathcal{J}(E,i).$\vspace*{2mm}

\noindent\textbf{Case 2.} Let $E=(f,g)_B,$ where $B$ is a cell in $\mathbb{R}^{d-1}$ and $f=-\infty$ or
$f:B\rightarrow\mathbb{R}$ is a continuous function with $[f,f]_B\subset\tilde{\pi}_d(W_{d-1}),$ and $g=+\infty$
or $g:B\rightarrow\mathbb{R}$ is a continuous function with $[g,g]_B\subset\tilde{\pi}_d(W_{d-1}).$ Moreover,
$(f,g)_B\cap\tilde{\pi}_d(W_{d-1})=\emptyset$ and either $f>-\infty$ or $g<+\infty$ (because
$\tilde{\pi}_d(W_{d-1})$ has surjective projection onto $\mathbb{R}^{d-1}$). Assume we already have constructed
$\mathcal{J}(B)=\{C_1,\ldots,C_s\}.$

For each $i\in\{1,\ldots,s\},$ let $N_i^1,\ldots,N_i^{j(i)}:\overline{C_i}\rightarrow\mathbb{R}$ be all continuous
maps with graphs contained in $\tilde{\pi}_d(W_{d-1}),$ whose restrictions to $B$ equal $g.$ If there are no such
maps, put $N_i=+\infty,$ otherwise, define $N_i:=\min\{N_i^1,\ldots, N_i^{j(i)}\}.$ Let
$M_i^1,\ldots,M_i^{k(i)}:\overline{C_i}\rightarrow\mathbb{R}$ be all continuous maps with graphs contained in
$\tilde{\pi}_d(W_{d-1}),$ whose restrictions to $B$ equal $f.$ If there are no such maps, put $M_i=-\infty,$
otherwise, define $M_i:=\max\{M_i^1,\ldots, M_i^{k(i)}\}.$

Define $\mathcal{J}(E)=\{(M_i,N_i)_{C_i}:i=1,\ldots,s\}.$\vspace*{2mm}
\begin{remark} The existence of $u,v$ in Case 1, follows by Proposition~\ref{basicprocell}(a)
\end{remark}

\begin{remark}Below we check that elements of $\mathcal{J}(E)$ are cells in $\mathbb{R}^d$ of dimension $d.$
By Proposition~\ref{basicprocell}(a), for $i=1,\ldots,s,$ $(C_i\times\mathbb{R})\cap\tilde{\pi}_d(W_{d-1})$ is the
union of pairwise disjoint graphs of continuous maps $C_i\rightarrow\mathbb{R},$ each of which has continuous
extension $\overline{C_i}\rightarrow\mathbb{R}$ (by Lemma \ref{closlem}(c)). The maps
$N_i^1,\ldots,N_i^{j(i)},M_i^1,\ldots,$$M_i^{k(i)}$ are among those extensions.

\end{remark}

\begin{remark}Induction on $d$ shows that the family $\mathcal{J}(E)$ is finite.
\end{remark}

\begin{lemma}\label{closeJ} Let $E$ be a cell in $\mathbb{R}^d.$ For each $C\in\mathcal{J}(E),$ $C$ is a cell
of dimension $d$ and $\overline{E}$ is a face of $\overline{C}.$
\end{lemma}
\noindent\textit{Proof.} Induction on $d.$ If $d=1,$ there is nothing to prove.

Let $d>1.$ First assume we are in Case 1 of the definition of $\mathcal{J}(E),$ that is $E=[f,f]_B.$ By induction
hypothesis, each element $C_i$ of $\mathcal{J}(B)$ is a cell. Let $T_{i,1}<\ldots<T_{i,t(i)}$ be all continuous
maps $C_i\rightarrow\mathbb{R}$ with graphs contained in $(M_i,N_i)_{C_i}\cap\tilde{\pi}_d(W_{d-1}).$ Put
$T_{i,0}:=M_i$ and $T_{i,t(i)+1}:=N_i.$ Clearly, each element of $\mathcal{J}(E)$ is of the form
$(T_{i,j},T_{i,j+1})_{C_i}$ for some $i=1,\ldots s$ and $j=0,\ldots,t(i),$ hence it is a cell of dimension
$\mathrm{dim}(C_i)+1$.

By induction hypothesis, $\overline{B}$ is a face of $\overline{C_i},$ for $i=1,\ldots,s.$ In particular,
$B\subset\overline{C_i}.$\vspace*{2mm}

\noindent\textbf{Claim a.} For $i=1,\ldots,s,$ we have:\vspace*{2mm}

\noindent(a) $v<+\infty$ if and only if $N_i<+\infty.$\\
(b) $u>-\infty$ if and only if $M_i>-\infty.$\vspace*{2mm}

\noindent In particular, $N_i|_B=v$ and $M_i|_B=u.$\vspace*{2mm}

\noindent\textit{Proof of Claim a.} (a) Fix any $i\in\{1,\ldots,s\}.$ If $N_i<+\infty,$ then $v<+\infty$ directly
by the definition of $N_i.$ Let $v<+\infty.$ Then, by $(3_{d-1})$ of Proposition~\ref{formulation},
Lemma~\ref{localsur} and Proposition~\ref{basicprocell}(a) and Lemma~\ref{closlem}(c), some of continuous maps
$C_i\rightarrow\mathbb{R}$ whose graphs are contained in $\tilde{\pi}_d(W_{d-1})$ extends to a continuous map
$\alpha:\overline{C_i}\rightarrow\mathbb{R}$ with $[v,v]_B\cap[\alpha,\alpha]_{\overline{C_i}}\neq\emptyset.$
Since $B$ is a cell, by Proposition~\ref{basicprocell}(a), $v=\alpha|_B$ which shows that
$\alpha\in\{N_i^1,\ldots,N_i^{j(i)}\}\neq\emptyset.$ Hence, $N_i<+\infty.$

The proof of (b) is very similar.\qed\vspace*{2mm}

By Lemma \ref{closlem}(c), each $T_{i,j}$ extends continuously to $\overline{C_i}\rightarrow\mathbb{R}$ also
denoted by $T_{i,j}$.\vspace*{2mm}

\noindent\textbf{Claim b.} For $i=1,\ldots,s$ and $j=1,\ldots,t(i),$ we have $T_{i,j}|_B=f|_B.$\vspace*{2mm}

\noindent\textit{Proof of Claim b.} Note that $$[T_{i,j},T_{i,j}]_B\subset [M_i,N_i]_{B}\cap
\tilde{\pi}_d(W_{d-1}),$$ so if $T_{i,j}|_B\neq f|_B,$ then, by Claim a, either $$v\leq T_{i,j}|_B\leq N_i|_B=v$$
or
$$u\geq T_{i,j}|_B\geq M_i|_B=u.$$ So either $T_{i,j}$ contradicts the choice of $N_i$ or the choice of
$M_i.$\qed\vspace*{2mm}

 Since, for $i=1,\ldots,s,$ $T_{i,0}|_B<f$ and $T_{i,t(i)+1}|_B>f,$ the list $T_{i,1},\ldots,T_{i,t(i)}$ is not
empty by $(3_{d-1})$ of Proposition~\ref{formulation}, Lemma~\ref{localsur} and Proposition~\ref{basicprocell}(a).
Combining Claim b with the induction hypothesis we get that $\overline{E}=[f,f]_{\overline{B}}$ is a face of
$[T_{i,j},T_{i,j}]_{\overline{C_i}}$, hence also a face of $[T_{i,j},T_{i,j+1}]_{\overline{C_i}},$ for
$j=0,\ldots,t(i),$ which completes the proof for $E$ as in Case 1. \vspace*{2mm}

Now assume we are in Case 2 of the definition of $\mathcal{J}(E),$ that is $E=(f,g)_B.$ By induction hypothesis,
each element $C_i$ of $\mathcal{J}(B)$ is a cell and $\overline{B}$ is a face of $\overline{C_i}$.\vspace*{2mm}

\noindent\textbf{Claim c.} For $i=1,\ldots,s,$ we have:\vspace*{2mm}

\noindent(a) $g<+\infty$ if and only if $N_i<+\infty,$\\
(b) $f>-\infty$ if and only if $M_i>-\infty.$\vspace*{2mm}

\noindent In particular, $N_i|_B=g$ and $M_i|_B=f.$\vspace*{2mm}

\noindent\textit{Proof of Claim c.}(a) Fix any $i\in\{1,\ldots,s\}.$ If $N_i<+\infty,$ then $g<+\infty$ directly
by the definition of $N_i.$

Assume $g<+\infty.$ Then, by $(3_{d-1})$ of Proposition~\ref{formulation}, Lemma~\ref{localsur} and
Proposition~\ref{basicprocell}(a) and Lemma~\ref{closlem}(c), some of continuous maps $C_i\rightarrow\mathbb{R}$
whose graphs are contained in $\tilde{\pi}_d(W_{d-1})$ extends to a continuous map
$\gamma:\overline{C_i}\rightarrow\mathbb{R}$ with $[g,g]_B\cap[\gamma,\gamma]_{\overline{C_i}}\neq\emptyset.$
Since $B$ is a cell, by Proposition \ref{basicprocell}(a), $g=\gamma|_B$ which shows that
$\gamma\in\{N_i^1,\ldots,N_i^{j(i)}\}\neq\emptyset.$ Hence, $N_i<+\infty.$

The proof of (b) is very similar.\qed\vspace*{2mm}

Now each element of $\mathcal{J}(E)$ is of the form $(M_i,N_i)_{C_i}$ for some $i=1,\ldots,s.$ First observe that
$(M_i,N_i)_{C_i}\cap\tilde{\pi}_d(W_{d-1})=\emptyset.$ Indeed, by induction hypothesis $C_i$ is a cell, so if
$(M_i,N_i)_{C_i}\cap\tilde{\pi}_d(W_{d-1})\neq\emptyset,$ then, by Proposition \ref{basicprocell}(a), there is a
continuous function $\beta: C_i\rightarrow\mathbb{R}$ whose graph is contained in $\tilde{\pi}_d(W_{d-1})$ with
$M_i<\beta<N_i.$ By Lemma \ref{closlem}(c), $\beta$ has continuous extension $\overline{C_i}\rightarrow\mathbb{R}$
denoted by the same letter. By induction hypothesis $B\subset\overline{C_i}.$ Now by Claim c either
$$g\leq\beta|_B\leq N_i|_B=g$$ or
$$f\geq\beta|_B\geq M_i|_B=f.$$ So $\beta$ contradicts the choice of $M_i$ or the choice of $N_i.$ We proved
$(M_i,N_i)_{C_i}\cap\tilde{\pi}_d(W_{d-1})=\emptyset.$ In particular, $(M_i,N_i)_{C_i}$ is a cell of dimension
$\mathrm{dim}(C_i)+1$.

Finally, Claim c implies that $[M_i|_B,N_i|_B]_B=[f,g]_B.$ By induction hypothesis, $\overline{B}$ is a face of
$\overline{C_i},$ so $\overline{E}=[f,g]_{\overline{B}}$ is a face of $[M_i,N_i]_{\overline{C_i}}.$\qed
\vspace*{3mm}

For a cell $E$ in $\mathbb{R}^d$ define $i(E):=\bigcup\{\overline{A}: A\in\mathcal{J}(E)\}.$

\begin{remark}\label{unionJ}Induction on $d$ shows that $$E\subset\mathrm{int}(i(E)).$$
\end{remark}

\begin{lemma} Let $E$ be a cell in $\mathbb{R}^d.$ Then  $\mathcal{J}(E)=\mathcal{I}(E).$
\end{lemma}
\noindent\textit{Proof.} By Lemma \ref{closeJ} we have $\mathcal{J}(E)\subset\mathcal{I}(E).$

Let $C$ be a cell of dimension $d$ which does not belong to $\mathcal{J}(E).$ Then for every $A\in\mathcal{J}(E)$
we have $\overline{A}\cap C=\emptyset.$ Indeed, by Lemmas \ref{frontier}, \ref{facclosedcell},  $\overline{A}$ is
the union of $A$ and cells of dimension smaller than $d.$ Then, by Proposition~\ref{basicprocell}(d),
$\overline{A}\cap C=\emptyset.$

The previous paragraph in combination with Remark~\ref{unionJ} gives $E\cap\overline{C}=\emptyset.$ Hence, $C$
does not belong to $\mathcal{I}(E).$ Thus $\mathcal{I}(E)\subset\mathcal{J}(E).$\qed

\begin{proposition}\label{isolation}Let $E$ be a cell in $\mathbb{R}^d$ of dimension $l.$ Let $C\in\mathcal{J}(E)$
and $F$ be a face of $\overline{C},$ $\mathrm{dim}(F)\leq l.$ If $F\neq \overline{E},$  then
$$F\cap\mathrm{int}(i(E))=\emptyset.$$
\end{proposition}
\textit{Proof.} By Lemmas \ref{facclosedcell}, \ref{frontier}, $F$ is the union of cells. If any $D$ of these
cells satisfies $D\cap\mathrm{int}(i(E))\neq\emptyset,$ then by Lemma \ref{boundandint} below,
$E\subset\overline{D}\subset F.$ Then $F=\overline{E}$ as $\mathrm{dim}(F)\leq\mathrm{dim}(E).$  \qed

\begin{lemma}\label{boundandint} Let $D, E$ be cells in $\mathbb{R}^d$ such that $D\cap\mathrm{int}(i(E))\neq\emptyset.$
Then $D\subset\mathrm{int}(i(E))$ and $E\subset\overline{D}.$
\end{lemma}
\textit{Proof.} Induction on $d.$ If $d=1,$ there is nothing to prove.

Assume that $d>1.$ Then either $E=[f,f]_B$ as in Case 1 or $E=(f,g)_B$ as in Case 2 of the construction of
$\mathcal{J}(E)$. Moreover, $A:=\pi_{d-1}(D)$ is a cell in $\mathbb{R}^{d-1}$ with
$A\cap\mathrm{int}(i(B))\neq\emptyset.$ By induction hypothesis, $A\subset\mathrm{int}(i(B))$ and
$B\subset\overline{A}.$

Assume first that $E=(f,g)_B$ (Case 2). According Case 2 of the construction of $\mathcal{J}(E)$,
$i(E)=\bigcup_{i=1}^s[M_i,N_i]_{\overline{C_i}}.$ Observe that
$$\hspace*{-30mm}(+)\hspace*{15mm}M_i|_{A}=M_j|_A\mbox{ and } N_i|_A=N_j|_A\mbox{ for all }i,j\mbox{ with
}A\subset\overline{C_i}\cap\overline{C_j}.$$ Indeed, by Claim c, $M_i|_B=M_j|_B.$ Suppose $M_i|_A\neq M_j|_A.$
Since $A$ is a cell we may assume that $M_i|_A<M_j|_A.$

As $A\subset\overline{C_i}$, by $(3_{d-1})$ of Proposition~\ref{formulation}, Lemma~\ref{localsur} and
Proposition~\ref{basicprocell}(a) and Lemma~\ref{closlem}(c), some of continuous maps $C_i\rightarrow\mathbb{R}$
whose graphs are contained in $\tilde{\pi}_d(W_{d-1})$ extends to a continuous map
$\kappa:\overline{C_i}\rightarrow\mathbb{R}$ with
$$[M_j|_A,M_j|_A]_A\cap[\kappa,\kappa]_{\overline{C_i}}\neq\emptyset.$$ Since $A$ is a cell, by Proposition
\ref{basicprocell}(a), $M_j|_A=\kappa|_A.$ But then $\kappa|_B=M_j|_B=M_i|_B$ and $\kappa|_{C_i}> M_i|_{C_i}.$
This contradicts the choice of $M_i$ in Case 2 of the construction of $\mathcal{J}(E).$ Hence $M_i|_A=M_j|_A$.
Similarly $N_i|_A=N_j|_A.$

Now $(+)$ in combination with the fact that $A\subset\mathrm{int}(i(A))\subset\mathrm{int}(i(B))$ implies that
$D=(M_j|_A,N_j|_A)_A$ (for any $j$ with $A\subset\overline{C_j}$) satisfies $D\subset\mathrm{int}(i(E)).$ Since
$E=(f,g)_B=(M_i|_B,N_i|_B)_B$ we have $E\subset \overline{D}.$

Assume next that $E=[f,f]_B$ as in Case 1 of the construction of $\mathcal{J}(E)$. Then
$i(E)=\bigcup_{i=1}^s[M_i,N_i]_{\overline{C_i}}$ (with $M_i, N_i$ different from those in Case 2). The same
argument as before shows that
$$\hspace*{-29mm}(-)\hspace*{15mm}M_i|_{A}=M_j|_A\mbox{ and } N_i|_A=N_j|_A\mbox{ for any }i,j\mbox{ with
}A\subset\overline{C_i}\cap \overline{C_j}.$$

Now $(-)$ in combination with the fact that $A\subset\mathrm{int}(i(A))\subset\mathrm{int}(i(B))$ implies that
$D\subset(M_j,N_j)_A$ (for any $j$ with $A\subset\overline{C_j}$) satisfies $D\subset\mathrm{int}(i(E)).$

By construction  in Case 1, $(M_j,N_j)_{C_j}\cap\tilde{\pi}_d(W_{d-1})$ is the union of graphs of functions on
$C_j$ having continuous extensions to $\overline{C_j}$ whose restrictions to $B$ equal $f.$ Therefore every cell
contained in $(M_j,N_j)_{A}$ contains $[f,f]_B$ in its closure.\qed

\begin{lemma}\label{boundunion}Let $E$ be a cell in $\mathbb{R}^d.$ Then $\mathrm{int}(i(E))$ and $\partial i(E)$ are the unions of cells.
\end{lemma}
\textit{Proof.} Since $i(E)$ is the union of closed cells, it is the union of cells by Lemma \ref{frontier}. Next
by Lemma \ref{boundandint}, $\mathrm{int}(i(E))$ is the union of cells. Therefore, by Proposition
\ref{basicprocell}(d), $\partial i(E)$ is the union of cells. \qed

\begin{proposition}\label{boundaryofiso}Let  $E$ be a cell in $\mathbb{R}^d.$ Then there is a polynomial
$P\in\mathbb{R}[x_1,\ldots,x_d]$ vanishing at every point of $\partial i(E),$ but at no point of $E.$
\end{proposition}
\textit{Proof.} Induction on $d.$ If $d=1,$ there is nothing to prove.

Let $d>1.$ Then $E=[f,f]_B$ as in Case 1 or  $E=(f,g)_B$ as in Case 2 of the construction of $\mathcal{J}(E).$
Here $B$ is a cell in $\mathbb{R}^{d-1}.$ We have $i(E)=\bigcup_{i=1}^s[M_i,N_i]_{\overline{C_i}},$ where $M_i,
N_i$ defined either as in Case 1 or as in Case 2 of the construction of $\mathcal{J}(E).$

Since $\mathrm{int}(i(B))$ is the union of cells (cf. Lemma \ref{boundunion}), by (+), (-) of the proof of Lemma
\ref{boundandint} we have that
$$M_i|_{\overline{C_i}\cap\overline{C_j}\cap\mathrm{int}(i(B))}=M_j|_{\overline{C_i}\cap\overline{C_j}\cap\mathrm{int}(i(B))}\mbox{ for }i, j=1,\ldots,s,$$
$$N_i|_{\overline{C_i}\cap\overline{C_j}\cap\mathrm{int}(i(B))}=N_j|_{\overline{C_i}\cap\overline{C_j}\cap\mathrm{int}(i(B))}\mbox{ for }i, j=1,\ldots,s.$$

It follows that
$$\partial i(E)\subset
(\partial
i(B)\times\mathbb{R})\cup\bigcup_{i=1}^s[M_i,M_i]_{\overline{C_i}}\cup\bigcup_{i=1}^s[N_i,N_i]_{\overline{C_i}}.$$

By induction hypothesis, there is a polynomial $Q\in\mathbb{R}[x_1,\ldots,x_{d-1}]$ vanishing at every point of
$\partial i(B),$ but at no point of $B,$ hence also at no point of $E.$ Recall that there is
$$q_{d-1}(x_1,\ldots,x_{d})=\prod_{l=0}^{\mathrm{deg}_{x_{d}}p_{d-1}}\frac{\partial^l p_{d-1}}{\partial x_{d}^l}(x_1,\ldots,x_{d}),$$
where $p_{d-1}\in(\mathbb{R}[x_1,\ldots,x_{d-1}])[x_{d}]$ is monic in $x_{d},$ such that, for $i=1,\ldots,s,$ we
have
$$[M_i,M_i]_{\overline{C_i}}, [N_i,N_i]_{\overline{C_i}}\subset\tilde{\pi}_{d}(W_{d-1})=\mathcal{Z}(q_{d-1}).$$
Hence, for $i=1,\ldots,s,$ there are $r(i), s(i)$ such that
$$[M_i,M_i]_{\overline{C_i}}\subset\mathcal{Z}(\frac{\partial^{r(i)} p_{d-1} }{\partial x_d^{r(i)}})\mbox{ and }
[N_i,N_i]_{\overline{C_i}}\subset\mathcal{Z}(\frac{\partial^{s(i)} p_{d-1} }{\partial x_d^{s(i)}}).$$

Observe that if $M_i>-\infty,$ then $\frac{\partial^{r(i)} p_{d-1} }{\partial x_d^{r(i)}}$ does not vanish at any
point of $E.$ Indeed, in Case 2 it follows immediately from the fact that $E\cap\tilde{\pi}_d(W_{d-1})=\emptyset.$
In Case 1 we have $E=[f,f]_B.$ Suppose that $$\frac{\partial^{r(i)} p_{d-1} }{\partial x_d^{r(i)}}(x,f(x))=0,
\mbox{ for some } x\in B.
$$ Since $$\frac{\partial^{r(i)} p_{d-1} }{\partial x_d^{r(i)}}(x,M_i(x))=0,$$ there is $M_i(x)<a<f(x)$ such that
$$\frac{\partial^{r(i)+1} p_{d-1} }{\partial x_d^{r(i)+1}}(x,a)=0.$$ In particular, the point $(x,a)$ satisfies
$$(x,a)\in (M_i|_B,f)_B\cap\tilde{\pi}_d(W_{d-1})=(u,f)_B\cap\tilde{\pi}_d(W_{d-1}).$$ This contradicts the fact
that $(u,f)_B\cap\tilde{\pi}_d(W_{d-1})=\emptyset.$

Similarly, if $N_i<+\infty,$ then $\frac{\partial^{s(i)} p_{d-1} }{\partial x_d^{s(i)}}$ does not vanish at any
point of $E.$

Finally, define $\check{P}_i:=\frac{\partial^{r(i)} p_{d-1} }{\partial x_d^{r(i)}},$ if $M_i>-\infty$ and
$\check{P}_i:=1,$ if $M_i=-\infty.$ Similarly, define $\hat{P}_i:=\frac{\partial^{s(i)} p_{d-1} }{\partial
x_d^{s(i)}},$ if $N_i<+\infty$ and $\hat{P}_i:=1,$ if $N_i=+\infty.$

Now it is sufficient to take $$P:=Q\prod_{i=1}^s\check{P}_i\cdot\hat{P}_i.$$The proof is complete.\qed

\begin{lemma}\label{zerooncell}Let $E\subset\mathbb{R}^d$ be a cell. Then there is
$Q\in\mathbb{R}[x_1,\ldots,x_d]$ such that $\mathcal{Z}(Q)\cap i(E)=\overline{E}.$
\end{lemma}
\textit{Proof.} Induction on $d.$ If $d=1,$ there is nothing to prove.

Let $d>1.$ Then $E=[f,f]_B$ as in Case 1 or  $E=(f,g)_B$ as in Case 2 of the construction of $\mathcal{J}(E).$
Here $B$ is a cell in $\mathbb{R}^{d-1}.$ By induction hypothesis there is $R\in\mathbb{R}[x_1,\ldots,x_{d-1}]$
such that $\mathcal{Z}(R)\cap i(B)=\overline{B}.$

In Case 2, by construction $\mathcal{Z}(R)\cap i(E)=\overline{E}.$

As for Case 1, recall that there is
$$q_{d-1}(x_1,\ldots,x_{d})=\prod_{l=0}^{\mathrm{deg}_{x_{d}}p_{d-1}}\frac{\partial^l p_{d-1}}{\partial x_{d}^l}(x_1,\ldots,x_{d}),$$
where $p_{d-1}\in(\mathbb{R}[x_1,\ldots,x_{d-1}])[x_{d}]$ is monic in $x_{d},$ such that,
$$E\subset\tilde{\pi}_{d}(W_{d-1})=\mathcal{Z}(q_{d-1}).$$ Hence, there is $l$ such that $E\subset\mathcal{Z}(\frac{\partial^l p_{d-1}}{\partial x_{d}^l}).$ As in the proof of
Proposition \ref{boundaryofiso} one shows that $\frac{\partial^l p_{d-1}}{\partial x_{d}^l}$ cannot vanish at any
point of $[u,u]_{\overline{B}}\setminus\overline{E}$ or $[v,v]_{\overline{B}}\setminus\overline{E}.$ Thus it is
sufficient to define $Q=R^2+(\frac{\partial^l p_{d-1}}{\partial x_{d}^l})^2.$\qed
\begin{remark}\label{zerocellrem} By Lemma \ref{boundandint}, Lemma \ref{frontier} and Remark \ref{unionJ}, $Q$ obtained in Lemma \ref{zerooncell}
satisfies $\mathcal{Z}(Q)\cap\mathrm{int}(i(E))=E.$
\end{remark}

\section{Proof of Theorem \ref{mainpw}}

\begin{proposition}\label{extentoiso}Let $E$ be a cell in $\mathbb{R}^n$ let $f:\overline{E}\rightarrow\mathbb{R}$ be a
continuous map such that $f|_E$ is regular, $f|_{\overline{E}\setminus E}=0.$ Then there exist
$P,Q\in\mathbb{R}[x_1,\ldots,x_n]$ such that $\mathcal{Z}(Q)\cap \mathrm{int}(i(E))=\emptyset,$
$\frac{P}{Q}|_E=f|_E,$ $\frac{P}{Q}:\mathrm{int}(i(E))\rightarrow\mathbb{R}$ extends continuously to
$i(E)\rightarrow\mathbb{R}$ by defining $\frac{P}{Q}|_{\partial i(E)}=0.$ Moreover,
$\mathrm{dim}\mathcal{Z}(Q)<\mathrm{dim}(E).$
\end{proposition}
\textit{Proof.} By Proposition \ref{boundaryofiso}, there is $w\in\mathbb{R}[x_1,\ldots,x_n]$ such that $$\partial
i(E)\subset\mathcal{Z}(w)\mbox{ and }\mathcal{Z}(w)\cap E=\emptyset.$$ By Lemma \ref{zerooncell} and Remark
\ref{zerocellrem} there is $v\in\mathbb{R}[x_1,\ldots,x_n]$ such that
$$\mathcal{Z}(v)\cap i(E)=\overline{E}\mbox{ and }\mathcal{Z}(v)\cap\mathrm{int}i(E)=E\mbox{ and } \mathcal{Z}(v)=\overline{E}^z.$$
It follows that $\overline{E}^z\cap i(E)=\overline{E}$ and $\overline{E}^z\cap\mathrm{int} i(E)={E}.$

Since $f$ is regular on $E,$ there are $p,q\in\mathbb{R}[x_1,\ldots,x_n]$ with $\mathcal{Z}(q)\cap E=\emptyset$
such that $f|_E=\frac{p}{q}|_E.$ The quotient $\frac{p}{q}$ defines a regular function on
$\overline{E}^z\setminus\mathcal{Z}(q)$ for which there is a regular function
$g:\mathbb{R}^n\setminus(\overline{E}^z\cap\mathcal{Z}(q))\rightarrow\mathbb{R}$ with
$g|_{\overline{E}^z\setminus\mathcal{Z}(q)}=\frac{p}{q}|_{\overline{E}^z\setminus\mathcal{Z}(q)}.$ In particular,
$g|_E=f|_E.$ Note that $\mathrm{int}i(E)\subset\mathbb{R}^n\setminus(\overline{E}^z\cap\mathcal{Z}(q)).$

\vspace*{2mm}

\noindent\textbf{Claim. }\textit{There exists a positive integer $m$ such that the regular function
$\tilde{g}:\mathrm{int}i(E)\rightarrow\mathbb{R}$ defined by $\tilde{g}=g\frac{w^{2m}}{w^{2m}+v^2}$
satisfies:}\vspace*{2mm}

\noindent(a) $\tilde{g}|_E=g|_E=f|_E$,\\
(b) $\tilde{g}$ \textit{ extends continuously to }$i(E)\rightarrow\mathbb{R}$ equal $0$ on $\partial
i(E).$\vspace*{2mm}

\noindent\textit{Proof of Claim.} (a) is obvious by the choice of $v, w.$

Let us turn to (b). Since $g$ is regular, there exist $s,t\in\mathbb{R}[x_1,\ldots,x_n]$ with
$\mathcal{Z}(t)\cap\mathbb{R}^n\setminus(\overline{E}^z\cap\mathcal{Z}(q))=\emptyset$ such that
$g=\frac{s}{t}|_{\mathbb{R}^n\setminus(\overline{E}^z\cap\mathcal{Z}(q))}.$ In particular, $$\mathcal{Z}(t)\cap
i(E) \subset i(E)\cap
(\overline{E}^z\cap\mathcal{Z}(q))=\overline{E}\cap\mathcal{Z}(q)\subset\overline{E}\setminus E.$$ Consequently,
by Theorem 2.6.6 of \cite{BCR}, there are a positive integer $k$ and a continuous semialgebraic function
$\alpha:i(E)\rightarrow\mathbb{R}$ such that $\alpha|_{\overline{E}\setminus E}=0$ and
$$\hspace*{-57mm}(\flat)\hspace*{55mm}\alpha(y)t(y)=d(y,\overline{E}\setminus E)^k$$ for $y\in i(E).$

Let $V$ be a semialgebraic tubular neighborhood of $E$ in $\mathrm{int}i(E).$ Let $\pi:V\rightarrow E$ be the
projection along the normal bundle. Let $d(y,E)$ denote the Euclidean distance from the point $y$ to $E.$ We may
assume that $|y-\pi(y)|=d(y,E),$ for $y\in V,$ where $|\cdot|$ is the Euclidean norm.

Define $$U=\{y\in V: d(y,E)<d(\pi(y),\overline{E}\setminus E)^{k}\}.$$ Our first aim is to check that\vspace*{2mm}

\noindent$(1)$  for every $c\in\overline{E}\setminus E$ and $\varepsilon>0$ there is $\delta>0$ such that for
every $y\in U$ with

$|y-c|<\delta$ we have $|g(y)|<\varepsilon.$\vspace*{2mm}

Clearly, there are positive constants $M, K$ such that for $y\in U$ with $|y-c|\leq 1$ we have
$$|s(y)|\leq |s(\pi(y))|+M|y-\pi(y)|=|s(\pi(y))|+Md(y,E)\leq$$
$$\leq|s(\pi(y))|+Md(\pi(y),\overline{E}\setminus E)^{k},$$

$$|t(y)|\geq |t(\pi(y))|-K|y-\pi(y)|=|t(\pi(y))|-Kd(y,E)\geq$$
$$\geq|t(\pi(y))|-Kd(\pi(y),\overline{E}\setminus E)^{k}\geq |t(\pi(y))|(1-K\alpha(\pi(y))),$$
the latter inequality by $(\flat).$ If $c\in\overline{E}\setminus E$ and $|y-c|<\delta,$ then
$$d(\pi(y),\overline{E}\setminus E)\leq |\pi(y)-c|\leq |y-c|+|y-\pi(y)|=|y-c|+d(y,E)\leq 2\delta.$$ Therefore,
if $\delta$ small enough, we have $|t(y)|\geq\frac{1}{2}t(\pi(y))$ and therefore
$$|g(y)|=|\frac{s(y)}{t(y)}|\leq 2|\frac{s(\pi(y))}{t(\pi(y))}|+2|\frac{Md(\pi(y),\overline{E}\setminus
E)^{k}}{t(\pi(y))}|= $$
$$=2|f(\pi(y))|+2M|\alpha(\pi(y))|.$$
Since $\pi(y)\in E,$ $f:\overline{E}\rightarrow\mathbb{R}$ and $\alpha|_{\overline{E}}$ are continuous,
$f(c)=\alpha(c)=0$ and $|\pi(y)-c|\leq 2\delta$ we obtain the condition $(1).$\vspace*{2mm}

Now we check the following:\vspace*{2mm}

\noindent (2) There exist a positive integer $l$ and a continuous semialgebraic function $\beta:i(E)\setminus
U\rightarrow\mathbb{R}$ such that for $m\geq l$ we have
$$\frac{w^{2m}}{w^{2m}+v^2}\leq w^{2(m-l)}\beta^2\mbox{ on }i(E)\setminus U.$$

First, note that on $i(E)\setminus U$ we  have $\mathcal{Z}(v)\subset\mathcal{Z}(w)$ hence by Theorem 2.6.6 of
\cite{BCR}, there are a positive integer $l$ and a continuous semialgebraic function $\beta:i(E)\setminus
U\rightarrow\mathbb{R}$ such that $w^l=\beta v$ on $i(E)\setminus U.$ Consequently, we have
$$\frac{w^{2m}}{w^{2m}+v^2}\leq \frac{\beta^2w^{2m}}{\beta^2w^{2m}+w^{2l}}\leq w^{2(m-l)}\beta^2,$$ as
required.\vspace*{2mm}

Let us complete the proof of (b). The extension of $\tilde{g}$ to $i(E)\rightarrow\mathbb{R}$ defined by assigning
$0$ to every point from $\partial i(E)$ will be denoted by $\hat{g}.$ Let us check that $\hat{g}$ is continuous at
every point of $\partial i(E)$. Fix $c\in \partial i(E).$ If $c\in\partial i(E)\setminus\overline{E},$ continuity
at $c$ is clear as $g$ is bounded in a neighborhood of $c$ and $v(c)\neq 0,$ $\frac{w^{2m}}{w^{2m}+v^2}(c)=0.$

Let $c\in\partial i(E)\cap\overline{E}\subset\overline{E}\setminus E.$ By (2)
$$\hat{g}\leq gw^{2(m-l)}\beta^2\mbox{ on }\mathrm{int}(i(E))\setminus U.$$
Since $\partial (i(E))\subset\mathcal{Z}(w),$ Proposition 2.6.4 of \cite{BCR} and the inequality above imply that
for $m$ large enough $\hat{g}|_{i(E)\setminus U}$ is continuous. Combining this with (1) and the fact that
$\frac{w^{2m}}{w^{2m}+v^2}$ is bounded we obtain that $\hat{g}$ is continuous at $c.$ The proof of Claim is
complete.\qed\vspace*{2mm}

Let us return to the proof of the proposition. Since $g$ is regular, there exist
$s,t\in\mathbb{R}[x_1,\ldots,x_n]$ with
$\mathcal{Z}(t)\cap\mathbb{R}^n\setminus(\overline{E}^z\cap\mathcal{Z}(q))=\emptyset$ such that
$g=\frac{s}{t}|_{\mathbb{R}^n\setminus(\overline{E}^z\cap\mathcal{Z}(q))}.$ Let $m$ be given by Claim. Put
$P=sw^{2m}$ and $Q=t(w^{2m}+v^2).$

In view of Claim, only the fact $\mathrm{dim}(\mathcal{Z}(Q))<\mathrm{dim}(E)$ requires explanation. The
inequality is an immediate consequence of $$\mathcal{Z}(t)\subset\overline{E}^z\cap\mathcal{Z}(q)$$ and
$$\mathcal{Z}(w^{2m}+v^2)=\mathcal{Z}(w)\cap\mathcal{Z}(v)=\mathcal{Z}(w)\cap\overline{E}^z$$ as $\overline{E}^z$
is an irreducible algebraic set and neither $q$ nor $w$ vanishes identically on $E.$ \qed\vspace*{3mm}

\noindent\textbf{Proof of Theorem \ref{mainpw}} Let $f:S\rightarrow\mathbb{R}$ be a piecewise-regular function.
Then there is a finite family $\{S_j\}_{j\in J}$ of semialgebraic subsets of $S$ such that $f|_{S_j}$ is regular
for each $j.$ After a linear change of coordinates in $\mathbb{R}^n,$ there exist real algebraic subsets
$W_0,\ldots,W_n$ of $\mathbb{R}^n$ satisfying the assertion of Proposition \ref{formulation}.

Now,  for each $d=1,\ldots,n,$ let $\mathcal{K}_d$ be a finite family  of cells in $\mathbb{R}^d$ as described in
Section \ref{constructofcells}. Let $\mathcal{K}=\mathcal{K}_{n}.$ By Propositions \ref{formulation},
\ref{basicprocell}(f), for every $j\in J$ and $C\in\mathcal{K}$ we have either $C\subset S_j$ or $C\cap
S_j=\emptyset.$

Let $\mathcal{K}^{(i)}=\{C\in\mathcal{K}: \mathrm{dim}(C)=i\}.$ Using induction on $i$ we will show that for
$0\leq i\leq \mathrm{dim}(S),$ there exists a continuous quasi-rational function
$g_i:\mathbb{R}^n\rightarrow\mathbb{R},$
$$g_i(x)=\frac{h_i(x)}{w_i(x)}, \mbox{ for  }x\in\mathbb{R}^n\setminus\mathcal{Z}(w_i),$$ where $h_i\in
QP(\mathbb{R}^n), w_i\in\mathbb{R}[x_1,\ldots,x_n]$ such that\vspace*{2mm}

(a) $g_i|_C${ is regular for each } $C\in\mathcal{K},$

(b) $g_i|_C=f|_C$  for each $C\in\bigcup_{l=0}^i\mathcal{K}^{(l)}$ with $C\subset S,$

(c) $w_i$ does not vanish identically on any $C\in\mathcal{K}$ with $\mathrm{dim}(C)\geq i.$\vspace*{2mm}

\noindent Once this is done, we complete the proof defining $g=g_{\mathrm{dim}(S)}.$

For $i=0,$ we can take $g_i$ to be any polynomial satisfying the interpolation problem (b). Assume that $i>0$ and
$g_{i-1}$ is already defined.

Let $v:S\rightarrow\mathbb{R}$ be given by the formula $v=f-g_{i-1.}$ Clearly, $v$ is a continuous function,
$v|_C$ is regular for every $C\in\mathcal{K}$ with $C\subset S,$ and $v|_C=0$ for every
$C\in\bigcup_{l=0}^{i-1}\mathcal{K}^{(l)}\mbox{ with }C\subset S.$

Pick any $E\in\mathcal{K}^{(i)}$ with $E\subset S.$ Then $\overline{E}\subset S$ because $S$ is closed. Moreover,
$v|_{\overline{E}\setminus E}=0$ by (b), because $\overline{E}\setminus E$ is the union of cells of dimension at
most $i-1$ as asserted by Lemmas \ref{frontier}, \ref{facclosedcell}.

Proposition \ref{extentoiso} applied to $v|_{\overline{E}}$ implies that there exists a continuous function
$u_E:\mathbb{R}^n\rightarrow\mathbb{R}$ such that $u_E|_{\mathrm{int}(i(E))}$ is regular and
$u_E|_{\mathbb{R}^n\setminus\mathrm{int} i(E)}=0$ and $u_E|_E=v|_E.$ More precisely, there exist $P_E,
Q_E\in\mathbb{R}[x_1,\ldots,x_n]$ with $\mathcal{Z}(Q_E)\cap\mathrm{int}(i(E))=\emptyset$  such that
$u_E|_{\mathrm{int}(i(E))}=\frac{P_E}{Q_E}|_{\mathrm{int}(i(E))}$ and
$\mathrm{dim}(\mathcal{Z}(Q_E))<\mathrm{dim}(E).$

Note that $u_E$ is regular on every cell in $\mathcal{K}.$ Indeed, by Lemma \ref{boundunion}, $\partial (i(E))$ is
the union of cells from $\mathcal{K}$. By Lemma \ref{frontier}, $i(E)$ is the union of cells from $\mathcal{K}$,
hence both $\mathrm{int}(i(E))$ and $\mathbb{R}^n\setminus\mathrm{int} i(E)$ are the unions of cells from
$\mathcal{K}.$ Now regularity of $u_E$ on every cell of $\mathcal{K}$ follows from the previous paragraph.

Finally, define $$g_i:=g_{i-1}+\sum_{E\in \mathcal{K}^{(i)}, E\subset S}u_E.$$ Clearly, $g_i|_C$ is regular for
every $C\in\mathcal{K}.$

To prove that $g_i$ satisfies (b), fix any $C\in\bigcup_{l=0}^{i}\mathcal{K}^{(l)}$ with $C\subset S.$ If
${C}\neq{E},$ then $C\cap\mathrm{int}(i(E))=\emptyset.$ Indeed, for $C\cap\mathrm{int}(i(E))\neq\emptyset$ we have
$E\subset\overline{C}$ by Lemma \ref{boundandint}. Then $E=C$ by Lemma \ref{frontier} and the fact that
$\mathrm{dim}(E)\geq\mathrm{dim}(C).$ Therefore, $u_E|_C=0$ for each $E\in\mathcal{K}^{(i)}$ with $E\subset S$ and
$E\neq C.$ Then if $\mathrm{dim}(C)\leq i-1,$ we have
$$g_i|_C=g_{i-1}|_C=f|_C.$$ If $\mathrm{dim}(C)=i,$ we have
$$g_i|_C=g_{i-1}|_C+u_C|_C=g_{i-1}|_C+v|_C=f|_C.$$ Hence (b) is satisfied.

Let us check that $g_i$ is a quotient of a quasi-polynomial and a polynomial function $w_i$ satisfying (c). It is
sufficient to show that each $u_E$ is such. Since $u_E|_{\partial i(E)}=0,$ we have $\partial
i(E)\subset\mathcal{Z}(P_E).$ Therefore $\tilde{P}_E$ defined to be $P_E$ on $i(E)$ and $-P_E$ on
$\mathbb{R}^n\setminus i(E)$ is a basic quasipolynomial function such that $$
u_E=\frac{P_E+\tilde{P}_E}{2Q_E}.$$\qed


\vspace*{3mm}

\noindent M. Bilski: Department of Mathematics and Computer Science,\\
Jagiellonian University,\\
\L ojasiewicza 6, 30-348 Krak\'ow, Poland.\\
E-mail: Marcin.Bilski@im.uj.edu.pl}\vspace*{2mm}\\

\end{document}